\documentclass[11pt,a4paper]{article}
\usepackage[utf8]{inputenc}
\usepackage[margin=2.6cm]{geometry}
\usepackage{amsmath,amssymb,amsthm,mathtools}
\usepackage{bm}
\usepackage{graphicx}
\usepackage[dvipsnames]{xcolor}
\usepackage{booktabs}
\usepackage{array}
\usepackage[section]{placeins}
\usepackage{enumitem}
\usepackage[numbers,sort&compress]{natbib}
\usepackage[colorlinks=true,linkcolor=blue!50!black,citecolor=green!40!black,urlcolor=blue!50!black]{hyperref}
\usepackage[capitalise,nameinlink]{cleveref}

\numberwithin{equation}{section}
\newtheorem{theorem}{Theorem}[section]
\newtheorem{proposition}[theorem]{Proposition}
\newtheorem{lemma}[theorem]{Lemma}

\theoremstyle{definition}
\newtheorem{definition}[theorem]{Definition}
\newtheorem{assumption}[theorem]{Assumption}
\theoremstyle{remark}
\newtheorem{remark}[theorem]{Remark}

\newcommand{\R}{\mathbb{R}}
\newcommand{\C}{\mathbb{C}}
\newcommand{\Z}{\mathbb{Z}}
\newcommand{\T}{\mathbb{T}}
\newcommand{\bu}{\bm{u}}
\newcommand{\bv}{\bm{v}}
\newcommand{\bw}{\bm{w}}
\newcommand{\be}{\bm{e}}
\newcommand{\bff}{\bm{f}}
\newcommand{\bz}{\bm{z}}
\newcommand{\ba}{\bm{a}}
\newcommand{\bk}{\bm{k}}
\newcommand{\bn}{\bm{n}}
\newcommand{\bphi}{\bm{\varphi}}
\newcommand{\bTheta}{\bm{\Theta}}
\newcommand{\dd}{\,\mathrm{d}}
\newcommand{\Vh}{\bm{V}_h}
\newcommand{\Qh}{Q_h}
\newcommand{\Ms}{\mathsf{m}}
\newcommand{\Wt}{\mathcal{W}}
\newcommand{\Rh}{\mathcal{R}_h}
\newcommand{\jump}[1]{[\![#1]\!]}
\newcommand{\norm}[1]{\|#1\|}
\newcommand{\snorm}[1]{|#1|}
\newcommand{\pih}{\pi_h}
\newcommand{\piF}{\pi_F}
\newcommand{\Ccal}{\mathcal{C}}
\newcommand{\E}{E_0}

\title{Assessing computational methods for large eddy simulation\\ through the spectrum of the residual}
\author{Erik Burman\thanks{Department of Mathematics, University College London, Gower Street, London WC1E 6BT, United Kingdom. \texttt{e.burman@ucl.ac.uk}.} \and Fabian Heimann\thanks{Department of Mathematics, University College London, Gower Street, London WC1E 6BT, United Kingdom. \texttt{f.heimann@ucl.ac.uk}}}
\date{October 26, 2026}

\begin{document}
\maketitle

\begin{abstract}
We propose a criterion for assessing computational methods for under-resolved simulations of incompressible flow (implicit large eddy simulation) that requires no regularity of the exact solution. In a dual weighted representation of the error in a linear quantity of interest the method enters only through its residual. Since the adjoint of the nonlinear problem couples all wave numbers, the unresolved part of the residual pollutes the large scales. For finite element methods on periodic domains, with velocities of polynomial degree $p$ on meshes of size $h$, let $\rho_{\bk}$ denote the Fourier coefficient of the residual for the wave vector $\bk$, tested in time against smooth weights. Using only energy stability, we prove bounds $\rho_{\bk}\lesssim h^{-\kappa}\min\{(h|\bk|)^{\sigma},1\}$, up to a viscous term and with constants independent of $h$, $\bk$ and the viscosity, so that resolved modes ($h|\bk|<1$) converge at the order $\sigma-\kappa$ and the unresolved coefficients are bounded by $h^{-\kappa}$, with tail exponent $\kappa\ge0$. Such bounds hold for Taylor--Hood elements with and without grad--div stabilisation, a nonconsistent streamline-diffusion method and the continuous interior penalty method with grad--div stabilisation; only the latter combines the resolved order $p+\frac12$ with the tail exponent $\frac12$. Computations for the Burgers equation and a double shear layer confirm the predicted behaviour for stable quantities of interest, and in forced two-dimensional turbulence the continuous interior penalty method, unlike the Taylor--Hood method with grad--div stabilisation alone, controls the grid-scale vorticity far from the viscous range.
\end{abstract}

\noindent\textbf{Keywords:} large eddy simulation; under-resolved turbulence; stabilised finite element methods; residual; Fourier analysis; dual weighted residual.

\smallskip
\noindent\textbf{MSC (2020):} 65M60, 65M15, 76F65, 76D05.

\section{Introduction}\label{sec:intro}

Resolving the Kolmogorov scale of a three-dimensional turbulent flow requires of the order of $Re^{9/4}$ grid points and a computational work of the order of $Re^{3}$, where $Re$ is the Reynolds number \cite{pope2000turbulent,moin1997tackling}, so that most flows of interest must be computed with only the largest scales resolved. In \emph{large eddy simulation} (LES) \cite{sagaut2006large} this is done either by filtering the equations and modelling the interaction between resolved and unresolved scales (physical LES), or by discretising the unmodified Navier--Stokes equations with a method that is robust enough to return a credible approximation on meshes that are coarse compared with the viscous scales (implicit LES, also called under-resolved direct numerical simulation (DNS)) \cite{grinstein2007implicit}. We are concerned with the second approach, and with the question of what, if anything, a given method approximates accurately before the viscous scales are resolved. This is the regime $\mu\ll Uh$, with viscosity $\mu$, a velocity scale $U$ and mesh size $h$, that is, a large mesh Reynolds number $Uh/\mu$; it distinguishes implicit LES from DNS. Our bounds hold uniformly in $\mu\ge0$, and in this regime their viscous contributions are smaller than the inviscid ones by a factor of order at most $(\mu/(Uh))^{1/2}$.

The lack of a stability theory for the Navier--Stokes equations in the inviscid limit has blocked theoretical progress on this question. In three space dimensions uniqueness of Leray--Hopf solutions is open, and non-uniqueness is known with forcing \cite{albritton2022nonuniqueness} and for weaker notions of solution \cite{buckmaster2019nonuniqueness}; weak solutions of the Euler equations are in general not unique \cite{delellis2009euler,delellis2010admissibility}, and even for smooth solutions the data-to-solution map is not uniformly continuous \cite{himonas2010nonuniform}. Error estimates in Sobolev norms rely on Gronwall's lemma, with constants $\exp(\int_0^T\norm{\nabla\bu}_{L^\infty}\dd t)$ for the exact velocity $\bu$ on the time interval $(0,T)$, and are void for turbulent flows. Convergence to suitable weak solutions \cite{guermond2005construction,guermond2006issues} does not discriminate between methods, since the standard Galerkin method has this property \cite{guermond2007faedo}, while it is considered unsuitable for under-resolved simulations. Existing criteria are based on dispersion--diffusion analysis of linear model problems \cite{moura2015linear,moura2020spatial,moura2024jointcg,moura2024jointdg}, on improved convergence orders for smooth solutions \cite{burman2006continuous,burman2007continuous,burman2015robust,defrutos2018analysis}, or on the requirement that the numerical dissipation should match the physical energy transfer \cite{burman2007interior,guasch2013statistical,fehn2023dissipation,fehn2025anomalous}; see, e.g., \cite{hughes2002variational,bazilevs2007weak,john2010numerical,schroeder2018divergence,fehn2018robust,lube2020implicit,fehn2022numerical} for finite element methods for turbulent flows. Hoffman and Johnson \cite{hoffman2004computability,hoffman2005computation,hoffman2006new,hoffman2007computational} used duality to argue that mean values may be computable when pointwise values are not, quantifying computability a posteriori through the size of an adjoint solution, in the spirit of the dual weighted residual method \cite{becker2001optimal,giles2002adjoint}.

We take the same duality argument as starting point, but use it to derive \emph{a priori} information on how the residual of a given method is distributed over wave numbers. We require of a criterion that it does not rely on regularity of the exact solution or on an asymptotic limit, that it applies to the nonlinear equations, that it can be tested theoretically and numerically, and that it distinguishes known poor methods from known good ones. Let $(\bu_h,p_h)$ be the discrete velocity and pressure, and let $(\cdot,\cdot)_\Omega$ denote the $L^2$ inner product on the spatial domain $\Omega$ and $\langle\cdot,\cdot\rangle$ a duality pairing. Assume that for a linear quantity of interest $J$ the error $\be=\bu-\bu_h$ satisfies
\begin{equation}\label{eq:intro_rep}
J(\be) = (\be(0),\bphi(0))_\Omega + \int_0^T \langle r_h(\bu_h,p_h)(t), (\bphi,\varsigma)(t)\rangle \dd t ,
\end{equation}
where $r_h$ is the residual of the discrete solution and $(\bphi,\varsigma)$, an adjoint velocity and pressure, solves a linearised adjoint problem (an additional consistency term appears when $\bu_h$ is not divergence free, \cref{prop:representation}). We assume that for reasonable methods the size of $\bphi$ is determined by the exact solution and $J$, so that the only method-specific quantity is $r_h$. For a linear problem with constant coefficients the error in the Fourier mode $\bk$ depends only on the residual in the same mode, and the large scales can be accurate even if the small-scale residual is large. For the Navier--Stokes equations the adjoint couples all wave numbers: on the torus $\Omega=\T^d:=\R^d/(2\pi\Z)^d$, $d\in\{2,3\}$, with the Fourier coefficients $\hat r_{\bk}$, $\hat\bphi_{\bk}$ for the wave vectors $\bk\in\Z^d$, $\langle r_h,\bphi\rangle=(2\pi)^{-d}\sum_{\bk}\hat r_{\bk}\cdot\overline{\hat\bphi_{\bk}}$, and slowly decaying high-frequency residual coefficients require a rapidly decaying adjoint spectrum, that is, a very stable quantity of interest.

\paragraph{Contributions.}
\begin{enumerate}[leftmargin=*]
\item We define the \emph{residual spectrum} $(\rho_{\bk})$ of a method through time-weak Fourier coefficients of its residual (\cref{def:spectrum}), prove \eqref{eq:intro_rep} (\cref{prop:representation}; see \cref{rem:formal} for its scope for turbulent flows), and show how bounds
\begin{equation}\label{eq:intro_bound}
\rho_{\bk} \le C h^{-\kappa} \min\{(h|\bk|)^{\sigma},1\}
\end{equation}
where $h$ is the mesh size and $C$ is independent of $h$, $\bk$ and the viscosity, combine with the decay of the adjoint spectrum to give convergence rates (\cref{thm:abstract,thm:negative}). The exponent $\sigma-\kappa$ is the convergence order of the resolved scales; $\kappa$ measures the amplification of the unresolved scales and sets the stability of $J$ required for convergence.
\item We prove bounds \eqref{eq:intro_bound}, uniform in the viscosity up to a viscous term and using only the energy estimate, for velocities of polynomial degree $p$: for the Taylor--Hood Galerkin method ($\kappa=1$, $\sigma-\kappa=p-1$), Taylor--Hood elements with grad--div stabilisation ($\kappa=1$, $\sigma-\kappa=p$), a nonconsistent streamline-diffusion method ($\kappa=\sigma-\kappa=\frac12$) and the continuous interior penalty (CIP) method with grad--div stabilisation ($\kappa=\frac12$, $\sigma-\kappa=p+\frac12$; for $p\ge2$ under an a posteriori verifiable condition), summarised in \cref{tab:summary}; $\ell^2$-type versions (\cref{def:l2type}) hold under a bound on $\norm{\bu_h}_{L^\infty}$.
\item Computations for the Burgers equation and for a double shear layer show the predicted residual spectra and the corresponding behaviour of the errors in the large scales when the quantity of interest is stable, and that for a sensitive quantity the error is small only once the mesh resolves the scales to which it is sensitive (\cref{sec:numerics}); in forced two-dimensional turbulence far from the viscous range the grad--div method loses control of the grid-scale vorticity under refinement, the CIP method does not (\cref{sec:forced}).
\end{enumerate}

\begin{remark}[Stable quantities]\label{rem:stable}
The bounds on the residual spectrum in \cref{sec:galerkin,sec:stabilised} are unconditional: they use only the energy estimate and hold uniformly in the viscosity, up to a common viscous term. What remains conditional is the convergence of the residual part of $J(\be)$, which by \cref{thm:abstract} requires the adjoint decay $\beta>d+\kappa$ ($s_*>\kappa$ in the norms of \cref{thm:negative}); this is a property of the flow and of $J$, not of the method (up to the dependence of the adjoint on $\bu_h$, \cref{rem:adjoint}), so the criterion compares methods for a given, possibly unknown, adjoint, and a smaller tail exponent $\kappa$ means that a less stable quantity still converges. No stable non-trivial quantities are known for high Reynolds number flow, so the argument cannot be closed in general: what is lacking is a stability estimate for a quantity of interest. When such an estimate is available, error estimates follow even for severely ill-posed problems \cite{burman2024data}; computationally, certain quantities appear to remain approximable as the viscosity vanishes \cite{hoffman2005computation,hoffman2007computational}, and time averages have been shown to be sufficiently stable to compare methods \cite{moura2022gradient}. For quantities of turbulent flows evaluated at a given time such estimates are in general available only with exponentially large constants (\cref{rem:formal}).
\end{remark}

The paper is organised as follows. \Cref{sec:model} treats one-dimensional model problems, \cref{sec:framework} the error representation, the residual spectrum and the abstract estimates, \cref{sec:galerkin,sec:stabilised} the four methods, \cref{sec:numerics} the numerical examples, and \cref{sec:conclusions} concludes.

\section{Model problems in one space dimension}\label{sec:model}

We first compare a linear convection--diffusion problem with the Burgers equation on the one-dimensional torus $\T=\R/2\pi\Z$. Throughout the paper, $\Omega$ denotes the spatial domain (here $\Omega=\T$), $(u,v)_\Omega := \int_\Omega u\,\bar v\dd x$ denotes the (complex) $L^2$ inner product, $\norm{\cdot}$ the corresponding norm, $\phi_k(x):=e^{ikx}$, $k\in\Z$, the Fourier modes and $\hat v_k := (v,\phi_k)_\Omega$ the (unnormalised) Fourier coefficients, so that $(u,v)_\Omega=(2\pi)^{-1}\sum_k \hat u_k\overline{\hat v_k}$. For a final time $T>0$ we let $I=(0,T)$, $Q=\T\times I$, $(u,v)_Q := \int_I (u,v)_\Omega\dd t$, and for $s\ge0$ and $x\ge0$ we write
\begin{equation}\label{eq:msdef}
\Ms_s(x) := \min\{x^s,1\}.
\end{equation}

Let $\{\mathcal T_h\}$ be a family of quasi-uniform meshes of $\T$ with mesh size $h\in(0,1]$ and let $V_h\subset H^1(\T)$ be the space of continuous, hence periodic, piecewise polynomials of degree $p\ge1$ on $\mathcal T_h$. We denote by $\pih$ the $L^2$-orthogonal projection onto $V_h$. We will use the following approximation properties of $\pih$ on Fourier modes, valid on quasi-uniform meshes thanks to the $L^\infty$-stability of $\pih$ \cite{douglas1974stability,crouzeix1987stability} (see \cref{lem:modes} below for the $d$-dimensional version); here and below $C$ denotes a generic constant independent of $h$, $k$ and of the diffusion coefficient $\mu$ introduced below:
\begin{equation}\label{eq:approx1d}
\norm{\phi_k-\pih\phi_k}_{L^\infty(\T)} \le C\,\Ms_{p+1}(h|k|),\qquad \norm{\partial_x(\phi_k-\pih\phi_k)}_{L^\infty(\T)} \le C\,|k|\,\Ms_{p}(h|k|).
\end{equation}

\subsection{A linear convection--diffusion equation}
Consider
\[
\partial_t u + \partial_x u - \mu\partial_{xx}u = f \quad\text{in } Q,\qquad u(\cdot,0)=u_0,
\]
with diffusion coefficient $\mu\ge0$, initial data $u_0\in L^2(\T)$ and source $f\in L^1(I;L^2(\T))$, and the Galerkin approximation: find $u_h\in W^{1,1}(I;V_h)$ with $u_h(0)=\pih u_0$ such that
\begin{equation}\label{eq:gal1d}
(\partial_t u_h + \partial_x u_h, v_h)_\Omega + \mu(\partial_x u_h,\partial_x v_h)_\Omega = (f,v_h)_\Omega\qquad \forall v_h\in V_h,\ \text{a.e. } t\in I.
\end{equation}
Testing with $v_h=u_h$ and using $(\partial_x u_h,u_h)_\Omega=0$ gives, with $\E := \norm{u_0}+\norm{f}_{L^1(I;L^2(\T))}$,
\begin{equation}\label{eq:energy1d}
\norm{u_h(t)}^2 + 2\mu\int_0^t\norm{\partial_x u_h}^2\dd s \le \E^2\qquad\forall t\in\bar I .
\end{equation}
(We detail the derivation of this type of bound in \cref{lem:energy}.) The constant is independent of $h$ and $\mu$, and hence so is the bound of the $L^2$-norm. This is essential for the discussion below, since the arguments must be independent of the Reynolds number to be relevant for LES.

We consider the Fourier coefficients of the error $e=u-u_h$, $\hat e_k(t)=(e(t),\phi_k)_\Omega$. To obtain an expression for $\hat e_k$ we introduce the adjoint problem
\begin{equation}\label{eq:adj1d}
-\partial_t\varphi - \partial_x\varphi - \mu\partial_{xx}\varphi = 0 \quad\text{in } Q,\qquad \varphi(\cdot,T)=\phi_k,
\end{equation}
whose solution is $\varphi(x,t) = e^{-\mu k^2(T-t)}\phi_k(x+T-t)$ (adjoint solutions are denoted by $\varphi$, Fourier modes by $\phi_k$). We define the residual of $u_h$ at time $t$ as the functional
\[
\langle r_h(t), v\rangle := (f - \partial_t u_h - \partial_x u_h, v)_\Omega - \mu(\partial_x u_h,\partial_x v)_\Omega,\qquad v\in H^1(\T),
\]
and its Fourier coefficients $\hat r_k(t) := \langle r_h(t),\phi_k\rangle$. Integrating by parts in space and time we obtain the exact representation
\begin{equation}\label{eq:rep1d}
\hat e_k(T) = e^{-(\mu k^2+ik)T}\,\hat e_k(0) + \int_0^T e^{-(\mu k^2+ik)(T-t)}\,\hat r_k(t)\dd t .
\end{equation}
The error at frequency $k$ depends only on the residual at the same frequency. By \eqref{eq:gal1d}, $\langle r_h(t),v_h\rangle=0$ for all $v_h\in V_h$, so that, with $w_k:=\phi_k-\pih\phi_k$ and since $(\partial_t u_h,w_k)_\Omega = 0$ and $(\partial_x u_h, w_k)_\Omega = -(u_h,\partial_x w_k)_\Omega$,
\[
\hat r_k(t) = \langle r_h(t),w_k\rangle = (f,w_k)_\Omega + (u_h,\partial_x w_k)_\Omega - \mu(\partial_x u_h,\partial_x w_k)_\Omega .
\]
By \eqref{eq:approx1d} and the Cauchy--Schwarz inequality,
\begin{equation}\label{eq:rk1d}
|\hat r_k(t)| \le C\big(\norm{f(t)}\,\Ms_{p+1}(h|k|) + (\norm{u_h(t)}+\mu\norm{\partial_x u_h(t)})\,|k|\,\Ms_p(h|k|)\big),
\end{equation}
and, combining with \eqref{eq:rep1d}, \eqref{eq:energy1d} and $|\hat e_k(0)|=|(u_0,w_k)_\Omega|\le C\norm{u_0}\Ms_{p+1}(h|k|)$,
\begin{equation}\label{eq:ek1d}
\sup_{t\in I}|\hat e_k(t)| \le C\,\big(1+T+(\mu T)^{1/2}\big)\,\E\;|k|\,\Ms_p(h|k|).
\end{equation}
Hence the low frequencies, that is the large scales, are approximated at the rate $O(h^p)$ without resolving the fine scales, independently of the viscosity and of the regularity of $u$; this holds in particular for the pure transport equation ($\mu=0$) with discontinuous data. For frequencies close to the resolution limit, $h|k|\sim1$, the bound grows like $|k|$: high-frequency content pollutes the fine-scale solution, but different frequencies do not interact, and the pollution is contained. For $h|k|\ge1$ the bound \eqref{eq:rk1d} reduces to the trivial bound $|\hat r_k(t)|\le C(\norm{f(t)}+(\norm{u_h(t)}+\mu\norm{\partial_xu_h(t)})|k|)$, showing that, after integration in time, the growth is at most linear in $|k|$ for the unresolved scales. For $\mu=0$ and smooth $f$ the coefficients themselves decay like $|k|^{-1}$ beyond the first aliasing wave numbers, since the residual is then piecewise polynomial with jumps at the nodes; the bound describes their envelope at the resolution limit (see the end of \cref{sec:burgers} and \cref{fig:burgers_spectra} for the Burgers equation).

\begin{remark}
The independence of $\mu$, and hence of the P\'eclet number, may be surprising, since spurious oscillations are expected to pollute the solution close to sharp layers in approximations using the standard Galerkin method. These spurious artefacts are still present and are expressed by the factor $|k|$ in \eqref{eq:ek1d}, but they remain at their own frequencies; as we shall see, this is not the case for nonlinear problems.
\end{remark}

\subsection{The Burgers equation}\label{sec:burgers}
We now replace the transport term by a nonlinear one,
\[
\partial_t u + \tfrac12\partial_x (u^2) - \mu\partial_{xx}u = f\quad\text{in }Q,\qquad u(\cdot,0)=u_0,
\]
and consider the Faedo--Galerkin method: find $u_h\in W^{1,1}(I;V_h)$, $u_h(0)=\pih u_0$, such that
\begin{equation}\label{eq:galburgers}
(\partial_t u_h + u_h\partial_x u_h, v_h)_\Omega + \mu(\partial_x u_h,\partial_x v_h)_\Omega = (f,v_h)_\Omega\qquad\forall v_h\in V_h,\ \text{a.e. } t\in I .
\end{equation}
Since $(u_h\partial_x u_h,u_h)_\Omega=\frac13\int_\T\partial_x(u_h^3)\dd x=0$, the bound \eqref{eq:energy1d} holds also in this case. To obtain a representation of $\hat e_k(T)$ we use the adjoint of the linearised (secant) error equation,
\begin{equation}\label{eq:adjburgers}
-\partial_t\varphi - u_{\rm m}\,\partial_x\varphi - \mu\partial_{xx}\varphi = 0\quad\text{in }Q,\qquad \varphi(\cdot,T)=\phi_k,\qquad u_{\rm m} := \tfrac12(u+u_h),
\end{equation}
which is exact since $\tfrac12(u^2-u_h^2)=u_{\rm m}\,e$. With $\langle r_h(t),v\rangle := (f-\partial_t u_h - u_h\partial_x u_h,v)_\Omega - \mu(\partial_x u_h,\partial_x v)_\Omega$ we arrive, as before, at
\[
\hat e_k(T) = (e(0),\varphi(0))_\Omega + \int_0^T\langle r_h(t),\varphi(t)\rangle\dd t .
\]
This time, however, $\varphi$ is no longer explicitly given and cannot be expected to be proportional to $\phi_k$, so that $\hat e_k$ does not depend only on $\hat r_k$. The reason becomes evident if we Fourier transform the adjoint equation \eqref{eq:adjburgers}: with $\hat\varphi_j(t) := (\varphi(t),\phi_j)_\Omega$,
\begin{equation}\label{eq:adjfourier}
-\partial_t\hat\varphi_j - \frac{1}{2\pi}\sum_{l\in\Z} (\hat u_{\rm m})_{j-l}\,(il)\,\hat\varphi_l + \mu j^2\hat\varphi_j = 0,\qquad \hat\varphi_j(T)=2\pi\delta_{jk} .
\end{equation}
If $u_{\rm m}$ has energy distributed over many frequencies, the convolution couples all wave numbers and amplifies the high-frequency content through the factor $il$. Expanding $\varphi(t)$ in its Fourier series we obtain
\begin{equation}\label{eq:repburgers}
\hat e_k(T) = (e(0),\varphi(0))_\Omega + \frac{1}{2\pi}\int_0^T\sum_{j\in\Z}\hat r_j(t)\,\overline{\hat\varphi_j(t)}\dd t ,
\end{equation}
and $\hat e_k(T)$ depends on all non-zero $\hat\varphi_j$, with a transfer between frequencies that depends on both the exact and the approximate solution. To conclude one would need bounds on $\hat\varphi_j$, uniform in $\mu$ and with precise decay in $j$, which are available only for simple model problems and methods designed for conservation laws \cite{tadmor1991local}. What we retain from \eqref{eq:repburgers} is that the decay of $\hat\varphi_j$ quantifies the difficulty of computing $\hat e_k$, while the only method-dependent contribution, which can be estimated rigorously, is the distribution of the residual over all frequencies $j$.

For the Galerkin method the residual coefficients are bounded as in the linear case: by \eqref{eq:galburgers}, $\hat r_k(t)=\langle r_h(t),w_k\rangle$ with $w_k=\phi_k-\pih\phi_k$, and since $|(u_h\partial_x u_h,w_k)_\Omega|\le \norm{u_h}\norm{\partial_xu_h}\norm{w_k}_{L^\infty}\le Ch^{-1}\norm{u_h}^2\norm{w_k}_{L^\infty}$ by an inverse inequality, we obtain
\begin{equation}\label{eq:rkburgers}
|\hat r_k(t)| \le C\Big(\norm{f(t)}\,\Ms_{p+1}(h|k|) + \norm{u_h(t)}^2\, h^{-1}\Ms_{p+1}(h|k|) + \mu\norm{\partial_x u_h(t)}\,|k|\,\Ms_p(h|k|)\Big).
\end{equation}
Note that $h^{-1}\Ms_{p+1}(h|k|) = |k|(h|k|)^p$ for $h|k|<1$ and $h^{-1}\Ms_{p+1}(h|k|)=h^{-1}\le|k|$ for $h|k|\ge1$. Splitting the sum in \eqref{eq:repburgers} at the frequency $N_c$ where the resolved-scale bound becomes $O(1)$, that is $N_c(hN_c)^p\sim1$, or $N_c\sim h^{-p/(p+1)}$, we obtain (with $\mu=0$ and $f=0$ for simplicity)
\begin{equation}\label{eq:splitburgers}
|\hat e_k(T) - (e(0),\varphi(0))_\Omega| \le C\int_0^T\Big(\underbrace{\sum_{|j|\le N_c}|j|(h|j|)^p|\hat\varphi_j(t)|}_{\text{resolved part}} + \underbrace{\sum_{|j|>N_c}\min\{|j|,h^{-1}\}\,|\hat\varphi_j(t)|}_{\text{unresolved part}}\Big)\dd t .
\end{equation}
The effect of the residual coefficients on the error differs from the linear case: if the $\hat\varphi_j$ are stable under refinement, the resolved part tends to zero, whereas the unresolved part, which represents the influence of the unresolved scales on $\hat e_k$, can only be made small by making $N_c$ large. Moving the derivative in $g_h:=u_h\partial_xu_h$ onto $\phi_j$ gives the trivial bound $\frac12\norm{u_h}^2|j|$, and leaving it gives $Ch^{-1}\norm{u_h}^2$. Writing $\phi_j=\partial_x(\phi_j/(ij))$ and integrating by parts on each element moves one more derivative onto $g_h$: $|(g_h,\phi_j)_\Omega|\le|j|^{-1}|Dg_h|(\T)$, where $|Dg_h|(\T)$ (jumps at the nodes plus $\norm{\partial_xg_h}_{L^1}$ on the elements) is at most $Ch^{-1}\norm{g_h}_{L^1}\le Ch^{-2}\norm{u_h}^2$ by inverse inequalities. Hence $|(g_h,\phi_j)_\Omega|\lesssim\norm{u_h}^2\min\{|j|,h^{-1},h^{-2}|j|^{-1}\}$, and with $\norm{\pih\phi_j}_{L^\infty}\lesssim(h|j|)^{-1}$ the same holds for $\hat r_j$ when $h|j|\ge1$. On a fixed mesh the residual spectrum thus has three regimes: growth like $|j|(h|j|)^p$ for $h|j|<1$, a peak of height $h^{-1}$ at $h|j|\approx1$, and decay like $(h|j|)^{-1}$; over a family of meshes the peaks lie on the envelope $|j|$, the factor in the unresolved part of \eqref{eq:splitburgers}. The rate of the decay is a property of the finite element space, not of the method, and the unresolved part remains dominated by $|j|\approx h^{-1}$. A method should therefore make $N_c$ as large as possible and reduce the amplification of the high modes, replacing the factor $|j|$ in the unresolved part by $|j|^\kappa$ with $\kappa<1$; a natural benchmark is the size of the projection error of a Fourier mode, $|\hat r_k|\lesssim\Ms_{p+1}(h|k|)$, that is $\kappa=0$ and resolved order $p+1$, which the CIP method attains in the example of \cref{sec:num1d}. For $\mu\lesssim Uh$, with $U$ a reference velocity, the viscous term in \eqref{eq:rkburgers} is bounded by $CU^{1/2}h^{1/2}|k|(h|k|)^p$ in $L^2(I)$, so that the convective term dominates.

\section{Navier--Stokes equations: error representation and residual spectrum}\label{sec:framework}

\subsection{Setting and notation}
Let $d\in\{2,3\}$, $\Omega=\T^d := \R^d/(2\pi\Z)^d$, $I=(0,T)$ with final time $T>0$, and $Q:=\Omega\times I$, with $(\cdot,\cdot)_Q$ and $\norm{\cdot}_Q$ the inner product and norm of $L^2(Q)$; for $m\ge0$, $|\cdot|_{W^{m,r}}$ and $|\cdot|_{H^m}$ denote the seminorms of the derivatives of order $m$. We use the notation of \cref{sec:model} also for vector-valued functions: $(\bu,\bv)_\Omega=\int_\Omega \bu\cdot\bar\bv\dd x$, $\phi_{\bk}(x)=e^{i\bk\cdot x}$ for $\bk\in\Z^d$, $\hat\bv_{\bk}=(\bv,\phi_{\bk})_\Omega\in\C^d$ and $(\bu,\bv)_\Omega=(2\pi)^{-d}\sum_{\bk}\hat\bu_{\bk}\cdot\overline{\hat\bv_{\bk}}$. $L^2_0(\Omega)$ denotes the functions in $L^2(\Omega)$ with zero mean, and $\mathrm{Id}$ the identity (matrix or operator). We consider the incompressible Navier--Stokes equations for the velocity $\bu$ and the pressure $p$ (in the rest of the paper $p$ denotes the polynomial degree, and $p_h$ the discrete pressure)
\begin{equation}\label{eq:ns}
\left\{\begin{aligned}
\partial_t\bu + (\bu\cdot\nabla)\bu + \nabla p - \mu\Delta\bu &= \bff &&\text{in } Q,\\
\nabla\cdot\bu &= 0 &&\text{in } Q,\\
\bu(\cdot,0) &= \bu_0 &&\text{in } \Omega,
\end{aligned}\right.
\end{equation}
with viscosity $\mu\ge 0$, body force $\bff\in L^1(I;[L^2(\Omega)]^d)$ and solenoidal initial velocity $\bu_0\in [L^2(\Omega)]^d$, and set $\E := \norm{\bu_0} + \norm{\bff}_{L^1(I;L^2(\Omega))}$. By a weak solution we mean a weakly solenoidal $\bu\in L^\infty(I;[L^2(\Omega)]^d)\cap C_w(\bar I;[L^2(\Omega)]^d)$ (where $C_w$ denotes weak continuity in time), with $\bu\in L^2(I;[H^1(\Omega)]^d)$ if $\mu>0$, such that, for all sufficiently smooth solenoidal test functions $\bv$,
\begin{equation}\label{eq:nsweak}
(\bu(T),\bv(T))_\Omega - (\bu_0,\bv(0))_\Omega - (\bu,\partial_t\bv)_Q - (\bu\otimes\bu,\nabla\bv)_Q + \mu(\nabla\bu,\nabla\bv)_Q = (\bff,\bv)_Q .
\end{equation}
For $\mu>0$ Leray--Hopf weak solutions exist; for $\mu=0$ weak solutions exist but are in general not unique \cite{delellis2010admissibility}. In two space dimensions and for $\mu>0$, weak solutions are unique and smooth for smooth data. In three space dimensions, uniqueness of Leray--Hopf solutions is open for $\bff=0$, and non-uniqueness has been shown for suitable forcing \cite{albritton2022nonuniqueness}; strong solutions are unique as long as they exist. In the inviscid limit continuous dependence on the data degenerates (see \cref{sec:intro}), so this case, which is the one of interest here, must be considered an ill-posed problem. In what follows $\bu$ denotes any weak solution; we never use its regularity or uniqueness.

For $s\ge0$ we use the function $\Ms_s$ of \eqref{eq:msdef}, and we measure temporal weights $\psi$ in
\[
\Wt := W^{1,1}(I;\R),\qquad \norm{\psi}_{\Wt} := \norm{\psi}_{L^\infty(I)} + \norm{\psi'}_{L^1(I)} .
\]
For complex or vector-valued functions $z$ of time, $\norm{z}_{\Wt}$ means $\norm{|z|}_{L^\infty(I)}+\norm{|z'|}_{L^1(I)}$.

\subsection{Adjoint representation of quantities of interest}
Let $\bu_1,\bu_2$ be solutions of \eqref{eq:ns} with data $(\bff_i,\bu_{0,i})$, $i=1,2$. The difference $\be=\bu_1-\bu_2$ satisfies a linear perturbation equation whose adjoint is, for given data $\bTheta_I\in L^1(I;[L^2(\Omega)]^d)$ and solenoidal $\bTheta_T\in[L^2(\Omega)]^d$: find an adjoint velocity and pressure $(\bphi,\varsigma)$ such that
\begin{equation}\label{eq:adjoint_general}
-\partial_t\bphi - (\bu_1\cdot\nabla)\bphi + (\nabla\bu_2)^{\top}\bphi - \nabla\varsigma - \mu\Delta\bphi = \bTheta_I,\qquad \nabla\cdot\bphi = 0,\qquad \bphi(T) = \bTheta_T .
\end{equation}
Whenever \eqref{eq:adjoint_general} has a sufficiently smooth solution, $(\be(T),\bTheta_T)_\Omega + (\be,\bTheta_I)_Q = (\be(0),\bphi(0))_\Omega + (\bff_1-\bff_2,\bphi)_Q$, so that the linear quantity $J(\bv):=(\bv(T),\bTheta_T)_\Omega+(\bv,\bTheta_I)_Q$ is stable, with a constant given by the size of $\bphi$ (and if this holds for a dense set of data, the solution of \eqref{eq:ns} is unique). The terminal data must be solenoidal since $\nabla\cdot\bphi=0$; a gradient part $\nabla\chi$, which does not affect $J$ on solenoidal fields, contributes $(\bu(T)-\bu_h(T),\nabla\chi)_\Omega=(\nabla\cdot\bu_h(T),\chi)_\Omega$ to the error of a discrete velocity $\bu_h$ of \cref{sec:fem}, a terminal mass term that is not part of the representation \eqref{eq:representation}. (The gradient part of $\bTheta_I$ is absorbed by $\varsigma$.)

\subsection{Finite element methods}\label{sec:fem}
Let $\{\mathcal T_h\}_{h\in(0,1]}$ be a family of conforming, shape-regular and quasi-uniform simplicial meshes of $\Omega$, compatible with periodicity, with mesh size $h$ (the maximal element diameter), and let $\mathcal F_h$ denote the set of (interior) faces. For $p\ge1$ let $V_h^p:=\{v\in H^1(\Omega):v|_K\in\mathbb P_p(K)\ \forall K\in\mathcal T_h\}$, where $\mathbb P_p(K)$ denotes the polynomials of degree at most $p$ on the element $K$, and consider
\[
\Vh := [V_h^p]^d,\qquad \Qh := V_h^{p_Q}\cap L^2_0(\Omega),\qquad p_Q\in\{p-1,p\},
\]
for the discrete velocity and pressure, with $p\ge2$ if $p_Q=p-1$ (Taylor--Hood) and $p_Q=p$ for the equal-order methods. Since the discrete velocity does not satisfy the divergence constraint exactly, the convective term is written in the skew-symmetric form
\begin{equation}\label{eq:skew}
c(\bw;\bu,\bv) := ((\bw\cdot\nabla)\bu,\bv)_\Omega + \tfrac12((\nabla\cdot\bw)\,\bu,\bv)_\Omega ,
\end{equation}
which satisfies $c(\bw;\bv,\bv)=0$ for real $\bw,\bv\in[H^1(\Omega)]^d$ and coincides with $((\bw\cdot\nabla)\bu,\bv)_\Omega=-(\bu\otimes\bw,\nabla\bv)_\Omega$ when $\nabla\cdot\bw=0$. All the methods considered in this paper are of the form: find $(\bu_h,p_h)\in W^{1,1}(I;\Vh)\times L^1(I;\Qh)$ with $\bu_h(0)=\bu_{0,h}$ such that, for almost every $t\in I$ and all $(\bv_h,q_h)\in\Vh\times\Qh$,
\begin{multline}\label{eq:fem}
(\partial_t\bu_h,\bv_h)_\Omega + c(\bu_h;\bu_h,\bv_h) - (p_h,\nabla\cdot\bv_h)_\Omega + (\nabla\cdot\bu_h,q_h)_\Omega + \mu(\nabla\bu_h,\nabla\bv_h)_\Omega \\
+ S_h(\bu_h;(\bu_h,p_h),(\bv_h,q_h)) = (\bff,\bv_h)_\Omega .
\end{multline}
Here $\bu_{0,h}\in\Vh$ is an $L^2$-stable approximation of $\bu_0$, $\norm{\bu_{0,h}}\le\norm{\bu_0}$, compatible with the mass equation, which contains no time derivative: $(\nabla\cdot\bu_{0,h},q_h)_\Omega=0$ for all $q_h$ that are not controlled by the stabilisation, that is, all $q_h\in\Qh$ for the Taylor--Hood methods and the $q_h$ in the kernel of the pressure stabilisation for the equal-order methods (for the CIP method of \cref{sec:sym} with $p\ge2$ this kernel contains the non-constant $C^1$ splines in $\Qh$). The discrete Leray projection used in \cref{sec:num2d} has these properties; the componentwise $L^2$-projection does not satisfy the constraint for Taylor--Hood elements. The stabilisation $S_h(\bw;\cdot,\cdot)$ is, for every $\bw$, a symmetric positive semi-definite bilinear form on $\Vh\times\Qh$; we write
\[
\snorm{(\bv_h,q_h)}_{S}^2:=S_h(\bu_h;(\bv_h,q_h),(\bv_h,q_h)).
\] All forms are extended sesquilinearly (conjugate linear in the last argument) to complex-valued arguments. The standard Galerkin method corresponds to $S_h=0$. We assume that $S_h$ depends smoothly on its first argument and that for the equal-order methods it controls the pressure up to a kernel of functions $q_h$ with $\nabla q_h\in\Vh$ (see \cref{sec:stabilised}), on which the inf-sup condition holds with $\bv_h=\nabla q_h$; then \eqref{eq:fem} is, after elimination of the pressure using the inf-sup condition and the pressure stabilisation, a system of ordinary differential equations (in the sense of Carath\'eodory) with locally Lipschitz right-hand side on the subspace of velocities satisfying the constraint above.

\begin{lemma}[Energy estimate]\label{lem:energy}
The solution of \eqref{eq:fem} exists for all $t\in\bar I$ and satisfies, for all $t\in\bar I$,
\begin{equation}\label{eq:energy}
\norm{\bu_h(t)}^2 + 2\int_0^t\big(\mu\norm{\nabla\bu_h}^2 + \snorm{(\bu_h,p_h)}_S^2\big)\dd s \le \E^2 .
\end{equation}
\end{lemma}
\begin{proof}
Taking $(\bv_h,q_h)=(\bu_h,p_h)$ in \eqref{eq:fem}, using $c(\bu_h;\bu_h,\bu_h)=0$ and the cancellation of the pressure terms, we obtain $\frac12\frac{\dd}{\dd t}\norm{\bu_h}^2 + \mu\norm{\nabla\bu_h}^2 + \snorm{(\bu_h,p_h)}_S^2 = (\bff,\bu_h)_\Omega\le\norm{\bff}\norm{\bu_h}$. Hence $\frac{\dd}{\dd t}(\norm{\bu_h}^2+\delta)^{1/2}\le\norm{\bff}$ for every $\delta>0$, and letting $\delta\to0$, $\norm{\bu_h(t)}\le M(t):=\norm{\bu_{0}}+\int_0^t\norm{\bff}\dd s$. Inserting this bound in the identity and integrating,
\[
\tfrac12\norm{\bu_h(t)}^2 + \int_0^t(\mu\norm{\nabla\bu_h}^2 + \snorm{(\bu_h,p_h)}_S^2)\dd s \le \tfrac12\norm{\bu_0}^2 + \int_0^t M'(s)M(s)\dd s = \tfrac12 M(t)^2\le\tfrac12\E^2 .
\]
Since all norms are equivalent on the finite dimensional space $\Vh$, the a priori bound implies that the local solution can be extended to all of $\bar I$.
\end{proof}
The constant in \eqref{eq:energy} is independent of $h$, $\mu$ and of the regularity of $\bu$; this is the only property of the discrete solution that will be used. Below, $C$ denotes a generic constant that may depend on $d$, $p$, the shape-regularity and quasi-uniformity of the meshes and the stabilisation parameters, but not on $h$, $\mu$, $T$, $\E$, $\bk$ or $\bu$. We write $\Ccal$ for a constant of the form $C(1+T)(1+\E)^2$ (for the stabilised methods depending also on the reference velocity $U$ of \cref{sec:stabilised}).

\subsection{Residual and residual spectrum}
For almost every $t\in I$ we define the residual of the discrete solution as the functional
\begin{multline}\label{eq:residual}
\langle r_h(t),(\bv,q)\rangle := (\bff,\bv)_\Omega - (\partial_t\bu_h,\bv)_\Omega - c(\bu_h;\bu_h,\bv) + (p_h,\nabla\cdot\bv)_\Omega - (\nabla\cdot\bu_h,q)_\Omega - \mu(\nabla\bu_h,\nabla\bv)_\Omega ,
\end{multline}
for $(\bv,q)\in[W^{1,\infty}(\Omega)]^d\times L^2(\Omega)$ (all arguments evaluated at time $t$). Note that $(\bu_h\cdot\nabla)\bu_h\in L^1(\Omega)$, so the convective term is well defined for $\bv\in L^\infty$. By construction, \eqref{eq:fem} is equivalent to the \emph{perturbed Galerkin orthogonality}
\begin{equation}\label{eq:orth}
\langle r_h(t),(\bv_h,q_h)\rangle = S_h(\bu_h;(\bu_h,p_h),(\bv_h,q_h))\qquad\forall(\bv_h,q_h)\in\Vh\times\Qh,\ \text{a.e. } t\in I .
\end{equation}
For a temporal weight $\psi$ and $(\bv,q)$ independent of time we set
\begin{equation}\label{eq:residual_tw}
\Rh(\psi;\bv,q) := \int_I \psi(t)\,\langle r_h(t),(\bv,q)\rangle\dd t .
\end{equation}
Since $\bu_h\in W^{1,1}(I;\Vh)$ we may integrate by parts in time in the term involving $\partial_t\bu_h$ when $\psi\in\Wt$; \cref{rem:timeweak} explains why this matters.

\begin{definition}[Residual spectrum]\label{def:spectrum}
For $\bk\in\Z^d$ let $\mathcal A_{\bk}:=\{\ba\in\C^d:\ \ba\cdot\bk=0,\ |\ba|=1\}$, so that $\ba\phi_{\bk}$ is solenoidal for $\ba\in\mathcal A_{\bk}$. The momentum and mass residual spectra are
\[
\rho^{\rm mom}_{\bk} := \sup_{\ba\in\mathcal A_{\bk}}\ \sup_{0\ne\psi\in\Wt}\frac{|\Rh(\psi;\ba\phi_{\bk},0)|}{\norm{\psi}_{\Wt}},\qquad
\rho^{\rm div}_{\bk} := \sup_{0\ne\psi\in L^2(I;\R)}\frac{|\Rh(\psi;\bm 0,\phi_{\bk})|}{\norm{\psi}_{L^2(I)}},
\]
and the residual spectrum is $\rho_{\bk}:=\rho^{\rm mom}_{\bk}+\rho^{\rm div}_{\bk}$.
\end{definition}

Since the momentum test functions are solenoidal, $\nabla\cdot(\ba\phi_{\bk})=i(\ba\cdot\bk)\phi_{\bk}=0$, the discrete pressure does not enter $\rho^{\rm mom}_{\bk}$. Only solenoidal test functions are needed in the momentum residual since the adjoint velocity is solenoidal; the mass residual is paired with the adjoint pressure, which requires no temporal regularity. We emphasise that $\rho_{\bk}$ is computable once a discrete solution is available (in practice for a finite set of temporal weights, see \cref{sec:numerics}).

\subsection{Error representation of quantities of interest}
Since $\bu_h$ is not divergence free, the linearisation of the convective term depends on the form in which it is written. For the representation of the error we use the conservative form of the exact equation, as in \eqref{eq:nsweak}. A direct computation using $\bu\otimes\bu-\bu_h\otimes\bu_h=\be\otimes\bu+\bu_h\otimes\be$, $\be=\bu-\bu_h$, shows that
\[
-(\bu\otimes\bu,\nabla\bv)_\Omega + (\bu_h\otimes\bu_h,\nabla\bv)_\Omega = \big(\be,\,-(\bu\cdot\nabla)\bv-(\nabla\bv)^{\top}\bu_h\big)_\Omega .
\]
The corresponding adjoint problem reads: find $(\bphi,\varsigma)$ such that
\begin{equation}\label{eq:adjoint}
-\partial_t\bphi - (\bu\cdot\nabla)\bphi - (\nabla\bphi)^{\top}\bu_h - \nabla\varsigma - \mu\Delta\bphi = \bTheta_I,\qquad \nabla\cdot\bphi = 0,\qquad \bphi(T) = \bTheta_T .
\end{equation}
Since $(\nabla\bu_h)^\top\bphi=\nabla(\bu_h\cdot\bphi)-(\nabla\bphi)^\top\bu_h$, \eqref{eq:adjoint} is of the form \eqref{eq:adjoint_general} with $\bu_1=\bu$, $\bu_2=\bu_h$, up to the gradient term $\nabla(\bu_h\cdot\bphi)$, which is absorbed in the multiplier. The discrete method, on the other hand, uses the skew-symmetric form; since $c(\bu_h;\bu_h,\bv) = -(\bu_h\otimes\bu_h,\nabla\bv)_\Omega - \tfrac12((\nabla\cdot\bu_h)\bu_h,\bv)_\Omega$, the difference gives rise to the \emph{divergence consistency term}
\begin{equation}\label{eq:Kh}
\epsilon_h(\bphi) := -\tfrac12\big((\nabla\cdot\bu_h)\,\bu_h,\bphi\big)_Q = \tfrac12\int_I\langle r_h(t),(\bm0,\bu_h\cdot\bphi)\rangle\dd t ,
\end{equation}
which vanishes for methods producing pointwise divergence-free velocities.

\begin{proposition}[Error representation]\label{prop:representation}
Let $\bu$ be a weak solution of \eqref{eq:ns}, $(\bu_h,p_h)$ the solution of \eqref{eq:fem} and $\be=\bu-\bu_h$. Assume that the adjoint problem \eqref{eq:adjoint} has a solution, satisfying \eqref{eq:adjoint} almost everywhere, with $\bphi\in W^{1,1}(I;[L^2(\Omega)]^d)\cap L^\infty(I;[W^{1,\infty}(\Omega)]^d)\cap L^2(I;[H^2(\Omega)]^d)$ and $\varsigma\in L^2(I;H^1(\Omega))$. Then, for $J(\bv)=(\bv(T),\bTheta_T)_\Omega+(\bv,\bTheta_I)_Q$,
\begin{equation}\label{eq:representation}
J(\be) = (\be(0),\bphi(0))_\Omega + \int_I\langle r_h(t),(\bphi(t),\varsigma(t))\rangle\dd t + \epsilon_h(\bphi).
\end{equation}
Moreover, with
\begin{equation}\label{eq:Phik}
\Phi_{\bk} := \norm{\hat\bphi_{\bk}}_{\Wt} + \norm{\hat\varsigma_{\bk}}_{L^2(I)},
\end{equation}
where $\hat\bphi_{\bk}(t)=(\bphi(t),\phi_{\bk})_\Omega$ and $\hat\varsigma_{\bk}(t)=(\varsigma(t),\phi_{\bk})_\Omega$, there holds
\begin{equation}\label{eq:spectral_bound}
\big|J(\be) - (\be(0),\bphi(0))_\Omega - \epsilon_h(\bphi)\big| \le 2d\,(2\pi)^{-d}\sum_{\bk\in\Z^d}\rho_{\bk}\,\Phi_{\bk}
\end{equation}
whenever the right-hand side is finite.
\end{proposition}
\begin{proof}
By mollification in space and time and dominated convergence (using $\bu\otimes\bu\in L^1(Q)$), \eqref{eq:nsweak} holds for all solenoidal $\bv$ in the class assumed for $\bphi$. We take $\bv=\bphi$ in \eqref{eq:nsweak}, and subtract the definition \eqref{eq:residual} tested with $(\bphi,\varsigma)$ and integrated over $I$, after an integration by parts in time of the term $(\partial_t\bu_h,\bphi)_\Omega$ and with $c(\bu_h;\bu_h,\bphi)$ written as above. Using $\nabla\cdot\bphi=0$ (so that the pressure terms vanish) and $-(\nabla\cdot\bu_h,\varsigma)_\Omega=(\nabla\cdot\be,\varsigma)_\Omega$, we obtain
\begin{multline*}
\int_I\langle r_h,(\bphi,\varsigma)\rangle\dd t = (\be(T),\bphi(T))_\Omega - (\be(0),\bphi(0))_\Omega - (\be,\partial_t\bphi)_Q - (\bu\otimes\bu-\bu_h\otimes\bu_h,\nabla\bphi)_Q \\ + \mu(\nabla\be,\nabla\bphi)_Q + (\nabla\cdot\be,\varsigma)_Q - \epsilon_h(\bphi).
\end{multline*}
By the identity preceding \eqref{eq:adjoint}, $\mu(\nabla\be,\nabla\bphi)_\Omega=(\be,-\mu\Delta\bphi)_\Omega$ and $(\nabla\cdot\be,\varsigma)_\Omega=(\be,-\nabla\varsigma)_\Omega$, the right-hand side equals $(\be(T),\bTheta_T)_\Omega-(\be(0),\bphi(0))_\Omega+(\be,\bTheta_I)_Q-\epsilon_h(\bphi)$, which is \eqref{eq:representation}. The second expression for $\epsilon_h$ in \eqref{eq:Kh} follows from $\langle r_h,(\bm0,q)\rangle=-(\nabla\cdot\bu_h,q)_\Omega$.

For \eqref{eq:spectral_bound}, let $S_N$ denote the Fourier partial sum over $|\bk|\le N$. For fixed $h$, the residual pairing in \eqref{eq:representation} is continuous with respect to solenoidal $\bphi$ in $L^\infty(I;L^2(\Omega))\cap L^1(I;H^1(\Omega))$ and to $\varsigma$ in $L^1(I;L^2(\Omega))$, since $\partial_t\bu_h\in L^1(I;L^2)$ and $(\bu_h\cdot\nabla)\bu_h,\nabla\bu_h\in L^\infty(Q)$, while the pressure term vanishes for solenoidal test functions; the partial sums $S_N\bphi$ are solenoidal. Since $\bphi\in C(\bar I;L^2(\Omega))$ and $\bphi(\bar I)$ is compact, $S_N\bphi\to\bphi$ uniformly in $t$ in $L^2(\Omega)$, and $S_N\bphi\to\bphi$ in $L^1(I;H^1(\Omega))$, $S_N\varsigma\to\varsigma$ in $L^2(I;L^2(\Omega))$ by dominated convergence. Hence the pairing is the limit of the partial sums $(2\pi)^{-d}\sum_{|\bk|\le N}\int_I\langle r_h,(\hat\bphi_{\bk}\phi_{\bk},\hat\varsigma_{\bk}\phi_{\bk})\rangle\dd t$. Since $\bk\cdot\hat\bphi_{\bk}=0$, let $\{\ba_{\bk,j}\}_{j}$ be a real orthonormal basis of $\bk^\perp$ ($d-1$ vectors if $\bk\neq0$, the canonical basis of $\R^d$ if $\bk=0$), and write $\hat\bphi_{\bk}(t)=\sum_j(\psi_{\bk,j}^{\rm re}(t)+i\psi_{\bk,j}^{\rm im}(t))\ba_{\bk,j}$ with real $\psi_{\bk,j}^{\rm re},\psi_{\bk,j}^{\rm im}\in\Wt$. By conjugate linearity,
\[
\Big|\int_I\langle r_h,(\hat\bphi_{\bk}\phi_{\bk},0)\rangle\dd t\Big| \le \sum_{j}\big(|\Rh(\psi^{\rm re}_{\bk,j};\ba_{\bk,j}\phi_{\bk},0)| + |\Rh(\psi^{\rm im}_{\bk,j};\ba_{\bk,j}\phi_{\bk},0)|\big)\le 2d\,\rho^{\rm mom}_{\bk}\,\norm{\hat\bphi_{\bk}}_{\Wt},
\]
and similarly $|\int_I\langle r_h,(\bm0,\hat\varsigma_{\bk}\phi_{\bk})\rangle\dd t|\le2\rho^{\rm div}_{\bk}\norm{\hat\varsigma_{\bk}}_{L^2(I)}$. This gives \eqref{eq:spectral_bound}.
\end{proof}


\begin{lemma}[Bounds for $\epsilon_h$]\label{lem:Kh}
For any method of the form \eqref{eq:fem},
\begin{equation}\label{eq:Kh_general}
|\epsilon_h(\bphi)|\le\tfrac12\,T^{1/2}\E\,\norm{\bphi}_{L^\infty(Q)}\,\norm{\nabla\cdot\bu_h}_{Q} .
\end{equation}
For the CIP method \eqref{eq:cip} of \cref{sec:sym}, with reference velocity $U$, semi-norms $\snorm{\cdot}_{S_u}$, $\snorm{\cdot}_{S_p}$, velocity and pressure dissipations $D_u:=\int_I\snorm{\bu_h}_{S_u}^2\dd t$, $D_p:=\int_I\snorm{p_h}_{S_p}^2\dd t$ and parameter $\varepsilon>0$ in \eqref{eq:cip_weight} (and $C$ depending on $\varepsilon^{-1/2}$),
\begin{equation}\label{eq:Kh_cip}
|\epsilon_h(\bphi)|\le C\Big(U^{-1/2}h^{1/2}\,T^{1/2}\E^2\,\norm{\bphi}_{L^\infty(I;W^{p+1,\infty}(\Omega))}+U^{-1}\norm{\bphi}_{L^\infty(Q)}\,(D_pD_u)^{1/2}\Big).
\end{equation}
\end{lemma}
\begin{proof}
The bound \eqref{eq:Kh_general} follows from the Cauchy--Schwarz inequality and \eqref{eq:energy}. For \eqref{eq:Kh_cip}, let $g:=\bu_h\cdot\bphi$, and see \cref{sec:galerkin,sec:sym}: $\pih$ is the $L^2$-projection onto $V_h^p$, $\jump{\cdot}$ the jump across a face and $\partial_{\bn}$ the derivative along a unit face normal $\bn$. By the mass equation, $(\nabla\cdot\bu_h,\pih g)_\Omega=-S_p(p_h,\pih g)$, so that
\[
2|\epsilon_h(\bphi)|\le\int_I\big(\norm{\nabla\cdot\bu_h}\norm{g-\pih g}+\snorm{p_h}_{S_p}\snorm{\pih g}_{S_p}\big)\dd t .
\]
On each element $K$, $g\in H^{p+1}(K)$, and inverse inequalities for $\bu_h$ give $\snorm{g}_{H^{p+1}(K)}\le C\sum_{j\le p}h^{-j}\norm{\bu_h}_{L^2(K)}\snorm{\bphi}_{W^{p+1-j,\infty}}$; hence the approximation properties of $\pih$ give the superapproximation estimate $\norm{g-\pih g}+h\norm{\nabla(g-\pih g)}\le Ch\norm{\bu_h}\norm{\bphi}_{W^{p+1,\infty}}$. Moreover $\jump{\partial_{\bn}g}=\jump{\partial_{\bn}\bu_h}\cdot\bphi$, hence $\jump{\partial_{\bn}\pih g}=\jump{\partial_{\bn}(\pih g-g)}+\jump{\partial_{\bn}\bu_h}\cdot\bphi$; using the trace inequality for the first term and the part $\varepsilon U$ of the weight in $S_u$ for the second (which introduces a factor $\varepsilon^{-1/2}$ in $C$),
\[
\snorm{\pih g}_{S_p}\le C\big(U^{-1/2}h^{1/2}\norm{\bu_h}\norm{\bphi}_{W^{p+1,\infty}}+U^{-1}\norm{\bphi}_{L^\infty}\snorm{\bu_h}_{S_u}\big).
\]
This, \cref{lem:div}, \eqref{eq:energy} and the bound $\int_I\snorm{p_h}_{S_p}\snorm{\bu_h}_{S_u}\dd t\le(D_pD_u)^{1/2}$ (Cauchy--Schwarz in time) give \eqref{eq:Kh_cip}.
\end{proof}

By \eqref{eq:Kh_general}, $|\epsilon_h|$ is bounded by $Ch^{-1}$ for the Taylor--Hood Galerkin method, by $C\gamma^{-1/2}$ for grad--div stabilisation with parameter $\gamma$ (\cref{thm:graddiv}), and by $Ch^{-1/2}$ for the stabilised methods of \cref{sec:stabilised}; for the CIP method with $\varepsilon>0$, \eqref{eq:Kh_cip} bounds it, up to a term of order $h^{1/2}$, by the geometric mean of the pressure and velocity dissipations, so that persistent velocity dissipation in under-resolved simulations does not by itself make it large. These are worst-case bounds and do not show that $\epsilon_h$ is large: since $\bphi$ is solenoidal, only the solenoidal part (Leray projection) of $(\nabla\cdot\bu_h)\bu_h$ enters, and its contribution depends on cancellations in space and time and on its pairing with the adjoint. In the estimates below we bound $J(\be)-(\be(0),\bphi(0))_\Omega-\epsilon_h(\bphi)$, assuming $\epsilon_h(\bphi)$ small. In \cref{sec:num2d} we compare the low-mode coefficients of $\frac12(\nabla\cdot\bu_h)\bu_h$, which enter $\epsilon_h$, with the bound \eqref{eq:Kh_general}.

\begin{remark}[On the adjoint]\label{rem:adjoint}
The proof of \cref{prop:representation} only inserts $\bphi$ in \eqref{eq:nsweak} and tests \eqref{eq:adjoint} with $\be$, and uses only the energy class of $\bu$. It therefore suffices that $\bphi\in W^{1,1}(I;L^2(\Omega))\cap L^1(I;W^{1,\infty}(\Omega))$ ($\cap\,L^2(I;H^1(\Omega))$ if $\mu>0$) and $\varsigma\in L^1(I;L^2(\Omega))$, with \eqref{eq:adjoint} holding weakly for all $\bv\in L^\infty(I;L^2(\Omega))$ ($\cap\,L^2(I;H^1(\Omega))$ if $\mu>0$) with $\nabla\cdot\bv\in L^\infty(I;L^2(\Omega))$, such as $\be$. In the transposed form \eqref{eq:adjoint} makes sense for every weakly solenoidal $\bu\in L^1(I;L^2(\Omega))$, but that form cannot be tested with $\be$; for $\mu=0$ the regularity above can be expected for general data only if the flow of $\bu$ is Lipschitz, $\nabla\bu\in L^1(I;L^\infty(\Omega))$, and this is where the regularity of $\bu$ enters. The estimates below need this regularity for each $h$, but only $\Phi_{\bk}$ (or $\norm{(\bphi,\varsigma)}_{X_s}$ in \cref{thm:negative}) enters the bounds. Since \eqref{eq:adjoint} is transported by $\bu$ and depends on $\bu_h$ through $(\nabla\bphi)^\top\bu_h$, its solution develops small scales even for smooth data, and its spectrum reflects those of $\bu$ and $\bu_h$.
Bounds on $\Phi_{\bk}$ uniform in $h$ and $\mu$ are therefore a stability assumption on the quantity of interest, which cannot be verified a priori for turbulent flows (\cref{rem:stable}). We assume in addition that the dependence of $\bphi$ on the discretisation is negligible for reasonable methods, so that the method-dependent part of the error is contained in $\rho_{\bk}$.
\end{remark}

\begin{remark}[Time-weak residuals]\label{rem:timeweak}
A pointwise-in-time version of \cref{def:spectrum} would require a bound on $\partial_t\bu_h$ itself, since $\bv-\piF\bv$, with $\piF$ the Fortin operator of \cref{ass:fortin}, is not orthogonal to $\Vh$, and \eqref{eq:energy} gives only $\norm{\partial_t\bu_h}\lesssim h^{-1-d/2}\E^2+\dots$. With time-weak residuals the temporal regularity of the adjoint enters instead, through $\Phi_{\bk}$. For the stabilised methods of \cref{sec:stabilised} the time derivative drops out exactly, and their bounds also hold for $\psi\in L^2(I)$ if $\bff\in L^2(I;L^2(\Omega))$.
\end{remark}

\subsection{An abstract estimate}\label{sec:abstract}
The following result shows how bounds on the residual spectrum and on the adjoint spectrum combine. We say that a method has \emph{residual spectrum of type $(\kappa,\sigma)$}, with exponents $0\le \kappa<\sigma$, if there are constants $C_R$ and $C_V$, independent of $h$, $\mu$ and $\bk$, such that for all $\bk\in\Z^d$ and $h\in(0,1]$
\begin{equation}\label{eq:type}
\rho_{\bk} \le C_R\,h^{-\kappa}\,\Ms_\sigma(h|\bk|) + C_V\,\mu^{1/2}\,|\bk|\,\Ms_p(h|\bk|) .
\end{equation}
The second term is the contribution of the viscous term; it is the same for all methods considered here and vanishes in the inviscid limit. For $h|\bk|<1$, $h^{-\kappa}\Ms_\sigma(h|\bk|)=h^{\sigma-\kappa}|\bk|^\sigma$, so $\sigma-\kappa$ is the convergence order of a fixed resolved mode, while for $h|\bk|\ge1$ the bound is the plateau $h^{-\kappa}$. Equivalently, since $h^{-\kappa}\Ms_\sigma(h|\bk|)\le|\bk|^\kappa\,\Ms_{\sigma-\kappa}(h|\bk|)$ for $|\bk|\ge1$, the unresolved coefficients are bounded uniformly in $h$ by the envelope $|\bk|^\kappa$.

\begin{theorem}\label{thm:abstract}
Assume that the method has residual spectrum of type $(\kappa,\sigma)$ and that the adjoint spectrum satisfies, as a stability assumption on the quantity of interest (\cref{rem:formal}), $\Phi_{\bk}\le C_A(1+|\bk|)^{-\beta}$ for some $\beta>d$, with $C_A$ independent of $h$ (the adjoint depends on $\bu_h$). Then, for $h\in(0,1/2]$,
\begin{equation}\label{eq:abstract}
\big|J(\be) - (\be(0),\bphi(0))_\Omega-\epsilon_h(\bphi)\big| \le C\,C_A\Big( C_R\,h^{-\kappa}\,\eta_{\sigma}(h) + C_V\,\mu^{1/2}\,h^{-1}\eta_{p+1}(h)\Big),
\end{equation}
where $\eta_s(h) := h^{s}$ if $\beta>d+s$, $\eta_s(h):=h^{s}|\log h|$ if $\beta=d+s$ and $\eta_s(h):=h^{\beta-d}$ if $d<\beta<d+s$; if $\mu>0$ we assume in addition $\beta>d+1$ (if $\mu=0$ the second term is absent). In particular, in the inviscid case the residual part of the error, $J(\be)-(\be(0),\bphi(0))_\Omega-\epsilon_h(\bphi)$, converges if $\beta>d+\kappa$, with rate $\min\{\sigma-\kappa,\beta-d-\kappa\}$ (up to a logarithm).
\end{theorem}
\begin{proof}
Let $N_n:=\#\{\bk\in\Z^d: n\le|\bk|<n+1\}\le C(1+n)^{d-1}$ and $n_h:=\lfloor 1/h\rfloor$. Splitting the sum in \eqref{eq:spectral_bound} into $|\bk|<1/h$ and $|\bk|\ge1/h$ we find
\begin{align*}
\sum_{|\bk|<1/h}h^{-\kappa}(h|\bk|)^\sigma(1+|\bk|)^{-\beta} &\le C h^{\sigma-\kappa}\sum_{n=0}^{n_h}(1+n)^{d-1+\sigma-\beta},\\
\sum_{|\bk|\ge1/h}h^{-\kappa}(1+|\bk|)^{-\beta}&\le Ch^{-\kappa}\sum_{n\ge n_h}(1+n)^{d-1-\beta}\le C h^{\beta-d-\kappa}.
\end{align*}
The first sum is bounded if $\beta>d+\sigma$, bounded by $C|\log h|$ if $\beta=d+\sigma$, and by $Ch^{\beta-d-\sigma}$ if $\beta<d+\sigma$. Since $\beta-d\ge\sigma$ in the first two cases, both contributions are bounded by $Ch^{-\kappa}\eta_\sigma(h)$. For the viscous term, $|\bk|\Ms_p(h|\bk|)=h^{p}|\bk|^{p+1}$ for $|\bk|<1/h$ and $=|\bk|$ otherwise, and the same computation gives the bound $C\big(h^p\sum_{n\le n_h}(1+n)^{d+p-\beta} + \sum_{n\ge n_h}(1+n)^{d-\beta}\big)\le Ch^{-1}\eta_{p+1}(h)$ for $\beta>d+1$.
\end{proof}

The theorem makes precise the heuristic discussion of \cref{sec:burgers}: the \emph{resolved order} $\sigma-\kappa$ determines the rate for very stable quantities ($\beta$ large), and the \emph{tail exponent} $\kappa$ the stability of the quantity of interest required for convergence and the rate for less stable ones. For a tolerance $0<\delta<C_Rh^{-\kappa}$, the modes for which the inviscid part of the bound in \eqref{eq:type} (its first term) is below $\delta$ are those with $|\bk|\le k_\delta(h):=(\delta/C_R)^{1/\sigma}h^{-\theta}$, with \emph{resolution exponent}
\begin{equation}\label{eq:theta}
\theta := \frac{\sigma-\kappa}{\sigma}\in(0,1].
\end{equation}
This is the cut-off frequency $N_c$ of \cref{sec:burgers}; the ideal value $\theta=1$ (all modes resolved by the mesh have small residual) is attained only if $\kappa=0$. In this sense a computational method for LES should have $\sigma-\kappa$ as large as possible and $\kappa$ as small as possible.

The per-mode hypothesis $\Phi_{\bk}\lesssim(1+|\bk|)^{-\beta}$ with $\beta>d+\kappa$ is, however, restrictive for turbulent flows: by \cref{rem:adjoint} the adjoint spectrum reflects the spectrum of the flow, and a velocity field with a Kolmogorov spectrum has Fourier coefficients of typical (shell-averaged) size $|\hat\bu_{\bk}|\sim|\bk|^{-11/6}$ when $d=3$ (and $|\bk|^{-2}$ in the enstrophy cascade range of two-dimensional turbulence), so that, if the adjoint inherits this spectrum, $\beta>d$ cannot be expected. The loss is due to the summation of per-mode bounds over shells containing $O(n^{d-1})$ modes. In the next subsection we show that, under a mild additional assumption on the discrete solution, the same exponents $(\kappa,\sigma)$ govern $\ell^2$-type (negative norm) bounds, for which the required adjoint regularity is reduced by $d/2$.

\subsection{Estimates in negative norms}\label{sec:negative}
For $s\ge0$ we use the Sobolev norms $\norm{v}_{H^s}^2:=(2\pi)^{-d}\sum_{\bk}(1+|\bk|^2)^s|\hat v_{\bk}|^2$, and for solenoidal $\bv$ and scalar $q$ defined on $Q$ we set
\[
\norm{(\bv,q)}_{X_s} := \norm{\bv}_{L^\infty(I;H^s)} + \norm{\partial_t\bv}_{L^1(I;H^s)} + \norm{q}_{L^2(I;H^s)} .
\]
We denote by $P_<$ and $P_\ge$ the Fourier projections onto the modes $|\bk|<1/h$ and $|\bk|\ge1/h$; they commute with $\partial_t$ and preserve solenoidality.

\begin{definition}\label{def:l2type}
Let $\mu=0$. A method has residual of \emph{$\ell^2$-type $(\kappa,\sigma)$} if there are constants $C_0,C_\sigma$, independent of $h$, such that for all $(\bv,q)$ with $\bv$ solenoidal (with the residual \eqref{eq:residual} extended to such test functions as in \cref{sec:negative_methods}),
\begin{equation}\label{eq:l2type}
\Big|\int_I\langle r_h,(\bv,q)\rangle\dd t\Big|\le \min\big\{C_0\,h^{-\kappa}\norm{(\bv,q)}_{X_0},\ C_\sigma\,h^{\sigma-\kappa}\norm{(\bv,q)}_{X_\sigma}\big\}.
\end{equation}
\end{definition}
The first bound is the $\ell^2$-analogue of the plateau $h^{-\kappa}$ for the unresolved modes, and the second of the bound $h^{\sigma-\kappa}|\bk|^\sigma$ for the resolved modes.

\begin{theorem}\label{thm:negative}
Let $\mu=0$ and assume that the method has residual of $\ell^2$-type $(\kappa,\sigma)$. Then, for all $0\le s\le\sigma$ and $h\in(0,1]$,
\[
\big|J(\be)-(\be(0),\bphi(0))_\Omega-\epsilon_h(\bphi)\big| \le C(C_0+C_\sigma)\,h^{s-\kappa}\,\norm{(\bphi,\varsigma)}_{X_s}.
\]
In particular, if $\norm{(\bphi,\varsigma)}_{X_s}$ is bounded uniformly in $h$ for all $s<s_*$ (a stability assumption, \cref{rem:formal}), the residual part of the error converges if $s_*>\kappa$, with rate $\min\{\sigma,s_*\}-\kappa-\delta$ for any $\delta>0$. If $\sup_t(|\hat\bphi_{\bk}(t)|+|\partial_t\hat\bphi_{\bk}(t)|)+\norm{\hat\varsigma_{\bk}}_{L^2(I)}\le C(1+|\bk|)^{-\beta}$ with $C$ independent of $h$, then $s_*=\beta-d/2$, so that $\beta>d/2+\kappa$ suffices, instead of $\beta>d+\kappa$ in \cref{thm:abstract}.
\end{theorem}
\begin{proof}
Write $(\bphi,\varsigma)=P_<(\bphi,\varsigma)+P_\ge(\bphi,\varsigma)$ in \eqref{eq:representation}. Since $1+|\bk|^2\ge h^{-2}$ for $|\bk|\ge1/h$ and $1+|\bk|^2\le 2h^{-2}$ for $|\bk|<1/h$, we have $\norm{P_\ge(\bphi,\varsigma)}_{X_0}\le h^s\norm{(\bphi,\varsigma)}_{X_s}$ and $\norm{P_<(\bphi,\varsigma)}_{X_\sigma}\le 2^{(\sigma-s)/2}h^{s-\sigma}\norm{(\bphi,\varsigma)}_{X_s}$. The first bound in \eqref{eq:l2type} applied to the high frequencies and the second to the low frequencies give the claim.
\end{proof}

In \cref{sec:negative_methods} we show that, if the discrete velocities are uniformly bounded, all the methods considered have residuals of $\ell^2$-type with the same exponents as their residual spectra. The gain has limits: for an adjoint whose spectrum mirrors a three-dimensional Kolmogorov spectrum, $s_*=\frac{11}{6}-\frac32=\frac13$, which is below $\kappa$ for all methods considered; in the enstrophy cascade range of two-dimensional turbulence $s_*=1$, which suffices for $\kappa=\frac12$ but not for the Galerkin method.

\begin{remark}[Stable quantities of interest for turbulent flows]\label{rem:formal}
The representation \eqref{eq:representation} is exact whenever \eqref{eq:adjoint} has a solution in the class of \cref{rem:adjoint}, and \cref{thm:abstract,thm:negative} hold under their hypotheses; for turbulent flows the question is whether these hold with useful constants. Some quantities at a given time are stable for every flow, for instance the mean velocity: for a constant vector $\ba$, $(\bu(T),\ba)_\Omega=(\bu_0,\ba)_\Omega+(\bff,\ba)_Q$. Quantities that couple to the unstable directions of a chaotic flow are not: their adjoint typically grows like $e^{\lambda(T-t)}$, with $\lambda$ the largest Lyapunov exponent, so that a perturbation of size $\delta$ affects them after a time of order $\lambda^{-1}\log(1/\delta)$ \cite{boffetta2002predictability}. In three dimensions $\lambda$ is at least of the order of the inverse of the Kolmogorov time scale, and therefore unbounded as $\mu\to0$ \cite{ruelle1979microscopic,boffetta2017chaos}; in this limit the time after which the large scales are affected stays bounded as $\delta\to0$ \cite{lorenz1969predictability}. For such quantities $C_A$ in \cref{thm:abstract} (and $\norm{(\bphi,\varsigma)}_{X_s}$ in \cref{thm:negative}) is of order $e^{\lambda T}$, although $|J(\be)|$ is bounded by the energy: the estimates remain valid for fixed $T$ and $\mu$, but are uninformative beyond the predictability horizon, and constants uniform in $\mu$ cannot be expected in three dimensions.

Averaged quantities, with $\bTheta_T=0$ and $\partial_t\bTheta_I=0$, may be stable with moderate constants: computations indicate that mean values may have dual solutions of moderate size \cite{hoffman2005computation,hoffman2007computational} and are shown to be robust enough to compare methods \cite{moura2022gradient}. By \cref{thm:negative,prop:negative}, a bound on $\norm{(\bphi,\varsigma)}_{X_s}$, uniform in $h$, for some $s>\frac12$ would suffice for the stabilised methods, and for some $s>1$ for the Taylor--Hood methods. For an example of stable averaged quantities in the context of a linear ill-posed problem, see \cite[Sect.~5.4]{burman2021hybridized}.
\end{remark}

\section{Galerkin methods with inf-sup stable elements}\label{sec:galerkin}

\subsection{Approximation of Fourier modes}
We first collect the approximation properties of discrete projections applied to Fourier modes. We denote by $\pih$ the $L^2$-orthogonal projection onto $V_h^p$ (applied componentwise to vector fields) and by $\pi_Q$ the $L^2$-projection onto $V_h^{p_Q}$. On quasi-uniform meshes $\pih$ is stable in $L^\infty(\Omega)$ \cite{douglas1974stability,crouzeix1987stability,diening2021sobolev}. Gradients of piecewise polynomials are understood elementwise, and for a face $F\in\mathcal F_h$, $\jump{\cdot}$ denotes the jump across $F$.

\begin{lemma}\label{lem:modes}
Let $\bk\in\Z^d$, $\ba\in\C^d$ with $|\ba|\le1$ and $\bv=\ba\phi_{\bk}$. Then
\begin{align}
\norm{\bv-\pih\bv}_{L^\infty(\Omega)} &\le C\,\Ms_{p+1}(h|\bk|), \label{eq:m1}\\
\norm{\nabla(\bv-\pih\bv)}_{L^\infty(\Omega)} &\le C\,|\bk|\,\Ms_{p}(h|\bk|), \label{eq:m2}\\
\norm{\nabla\pih\bv}_{L^\infty(\Omega)} &\le C\,h^{-1}\Ms_{1}(h|\bk|), \label{eq:m3}\\
\max_{F\in\mathcal F_h}\norm{\jump{\nabla\pih\bv}}_{L^\infty(F)} &\le C\,h^{-1}\Ms_{p+1}(h|\bk|), \label{eq:m4}
\end{align}
and the same bounds hold for $\pi_Q\phi_{\bk}$ with $p$ replaced by $p_Q$. The corresponding $L^2$-norms are bounded by $(2\pi)^{d/2}$ times the right-hand sides.
\end{lemma}
\begin{proof}
All derivatives of order $m$ of $\bv$ are bounded by $|\bk|^m$. If $h|\bk|<1$, the $L^\infty$-stability of $\pih$ and the approximation properties of the Lagrange interpolant $I_h$ give $\norm{\bv-\pih\bv}_{L^\infty}\le C\inf_{\bv_h\in\Vh}\norm{\bv-\bv_h}_{L^\infty}\le C(h|\bk|)^{p+1}$, and by an inverse inequality
\[
\norm{\nabla(\bv-\pih\bv)}_{L^\infty}\le\norm{\nabla(\bv-I_h\bv)}_{L^\infty}+Ch^{-1}\big(\norm{I_h\bv-\bv}_{L^\infty}+\norm{\bv-\pih\bv}_{L^\infty}\big)\le C|\bk|(h|\bk|)^p .
\]
If $h|\bk|\ge1$, the $L^\infty$-stability gives $\norm{\bv-\pih\bv}_{L^\infty}\le C$ and, by an inverse inequality, $\norm{\nabla\pih\bv}_{L^\infty}\le Ch^{-1}\norm{\pih\bv}_{L^\infty}\le Ch^{-1}\le C|\bk|$, so that $\norm{\nabla(\bv-\pih\bv)}_{L^\infty}\le C|\bk|$. This proves \eqref{eq:m1}--\eqref{eq:m2}. For \eqref{eq:m3} we use $\norm{\nabla\pih\bv}_{L^\infty}\le\norm{\nabla\bv}_{L^\infty}+\norm{\nabla(\bv-\pih\bv)}_{L^\infty}\le C|\bk|$ together with $\norm{\nabla\pih\bv}_{L^\infty}\le Ch^{-1}$. Finally, since $\bv$ is smooth, $\jump{\nabla\pih\bv}=\jump{\nabla(\pih\bv-\bv)}$, which is bounded by $2\norm{\nabla(\bv-\pih\bv)}_{L^\infty}\le C|\bk|(h|\bk|)^p=Ch^{-1}(h|\bk|)^{p+1}$ for $h|\bk|<1$ and by $2\norm{\nabla\pih\bv}_{L^\infty}\le Ch^{-1}$ in general.
\end{proof}

For the Taylor--Hood element we use a Fortin operator.
\begin{assumption}[Fortin operator]\label{ass:fortin}
There is a linear operator $\piF:[W^{1,1}(\Omega)]^d\to\Vh$ such that $(q_h,\nabla\cdot(\bv-\piF\bv))_\Omega=0$ for all $q_h\in\Qh$ and, for $1\le m\le p+1$,
\[
\norm{\bv-\piF\bv}_{L^\infty(\Omega)}+h\norm{\nabla(\bv-\piF\bv)}_{L^\infty(\Omega)}\le C_F\,h^m\,|\bv|_{W^{m,\infty}(\Omega)} .
\]
\end{assumption}
For the lowest-order Taylor--Hood element ($p=2$) on shape-regular simplicial meshes in any dimension, such an operator, which is moreover locally $W^{1,r}$-stable for all $r\in[1,\infty]$, was constructed in \cite{diening2022fortin} (the operator is a local projection onto $\Vh$, built from a Scott--Zhang operator and a divergence correction; the construction applies verbatim to periodic meshes, on which every simplex has an interior vertex). For $\bv=\ba\phi_{\bk}$ and $h|\bk|<1$, \cref{ass:fortin} with $m=p+1$ gives
\begin{equation}\label{eq:fortin_modes}
\norm{\bv-\piF\bv}_{L^\infty(\Omega)}\le C(h|\bk|)^{p+1},\qquad \norm{\nabla(\bv-\piF\bv)}_{L^\infty(\Omega)}\le C|\bk|(h|\bk|)^{p}.
\end{equation}

\subsection{The Taylor--Hood Galerkin method}
We consider \eqref{eq:fem} with $p\ge2$, $p_Q=p-1$ and $S_h=0$, the Taylor--Hood (TH) Galerkin method. For $p\ge3$ \cref{ass:fortin} is an assumption (the inf-sup condition alone only provides an $H^1$-stable, in general nonlocal, Fortin operator). We denote by
\[
\mathcal V_{\bk}:=C(\mu T)^{1/2}\E\,|\bk|\,\Ms_p(h|\bk|)
\]
the viscous contribution, the second term in \eqref{eq:type}.

\begin{theorem}\label{thm:th}
Let \cref{ass:fortin} hold. The residual spectrum of the Taylor--Hood Galerkin method satisfies, for all $\bk\in\Z^d$, $h\in(0,1]$ and $\mu\ge0$,
\begin{align}
\rho^{\rm mom}_{\bk} &\le \Ccal\,h^{-1}\Ms_{p+1}(h|\bk|) + \mathcal V_{\bk}, \label{eq:th_mom}\\
\rho^{\rm div}_{\bk} &\le C\,T^{1/2}\E\,h^{-1}\Ms_{p}(h|\bk|). \label{eq:th_div}
\end{align}
In particular the method has residual spectrum of type $(\kappa,\sigma)=(1,p)$, with resolved order $\sigma-\kappa=p-1$; the momentum residual alone is of type $(1,p+1)$.
\end{theorem}
\begin{proof}
Let $\ba\in\mathcal A_{\bk}$, $\bv=\ba\phi_{\bk}$ and $\psi\in\Wt$. We use repeatedly that, by \eqref{eq:energy}, $\sup_t\norm{\bu_h(t)}\le\E$ and $\mu\norm{\nabla\bu_h}_Q^2\le\E^2/2$, and the inverse inequality $\norm{\nabla\bu_h}\le C_Ih^{-1}\norm{\bu_h}$. For any time-independent $\bw\in[W^{1,\infty}(\Omega)]^d$ with $(p_h,\nabla\cdot\bw)_\Omega=0$ we have
\begin{multline*}
\int_I\psi\langle r_h,(\bw,0)\rangle\dd t = \underbrace{\int_I\psi(\bff,\bw)_\Omega\dd t}_{\mathsf T_1}\ \underbrace{-\int_I\psi(\partial_t\bu_h,\bw)_\Omega\dd t}_{\mathsf T_2}\\ \underbrace{-\int_I\psi\, c(\bu_h;\bu_h,\bw)\dd t}_{\mathsf T_3}\ \underbrace{-\mu\int_I\psi(\nabla\bu_h,\nabla\bw)_\Omega\dd t}_{\mathsf T_4},
\end{multline*}
and the four terms are bounded as follows. Clearly $|\mathsf T_1|\le\norm{\psi}_{L^\infty}\norm{\bff}_{L^1(I;L^2)}\norm{\bw}$. Integrating by parts in time,
\[
|\mathsf T_2| = \Big|\psi(T)(\bu_h(T),\bw)_\Omega-\psi(0)(\bu_h(0),\bw)_\Omega-\int_I\psi'(\bu_h,\bw)_\Omega\dd t\Big|\le 2\norm{\psi}_{\Wt}\E\norm{\bw}.
\]
Since $|(\bu_h\cdot\nabla)\bu_h|+\frac12|(\nabla\cdot\bu_h)\bu_h|\le(1+\frac{\sqrt d}{2})|\bu_h||\nabla\bu_h|$ pointwise,
\[
|c(\bu_h;\bu_h,\bw)|\le C\norm{\bu_h}\norm{\nabla\bu_h}\norm{\bw}_{L^\infty}\le Ch^{-1}\norm{\bu_h}^2\norm{\bw}_{L^\infty},
\]
so that $|\mathsf T_3|\le C\,T\,\norm{\psi}_{L^\infty}\E^2h^{-1}\norm{\bw}_{L^\infty}$.
Finally $|\mathsf T_4|\le\norm{\psi}_{L^\infty}\mu^{1/2}T^{1/2}(\mu\norm{\nabla\bu_h}_Q^2)^{1/2}\norm{\nabla\bw}\le\norm{\psi}_{L^\infty}(\mu T)^{1/2}\E\norm{\nabla\bw}$. Altogether, using $h\le1$,
\begin{equation}\label{eq:basic}
\Big|\int_I\psi\langle r_h,(\bw,0)\rangle\dd t\Big|\le C\norm{\psi}_{\Wt}\Big((1+T)(\E+\E^2)\,h^{-1}\norm{\bw}_{L^\infty}+(\mu T)^{1/2}\E\,\norm{\nabla\bw}\Big).
\end{equation}

\emph{Case $h|\bk|\ge1$.} Since $\nabla\cdot\bv=0$ we may take $\bw=\bv$ in \eqref{eq:basic}. With $\norm{\bv}_{L^\infty}=|\ba|=1$ and $\norm{\nabla\bv}\le(2\pi)^{d/2}|\bk|$ this gives $|\Rh(\psi;\bv,0)|\le\norm{\psi}_{\Wt}(\Ccal h^{-1}+\mathcal V_{\bk})$, which is \eqref{eq:th_mom} in this case. Note that neither the Galerkin orthogonality nor the Fortin operator is used here: the bound only expresses the $L^2$-stability of the method.

\emph{Case $h|\bk|<1$.} By \eqref{eq:orth} with $S_h=0$ we have $\langle r_h(t),(\piF\bv,0)\rangle=0$, hence $\Rh(\psi;\bv,0)=\int_I\psi\langle r_h,(\bw,0)\rangle\dd t$ with $\bw:=\bv-\piF\bv$, and $(p_h,\nabla\cdot\bw)_\Omega=0$ by the Fortin property since $p_h\in\Qh$. Inserting \eqref{eq:fortin_modes} into \eqref{eq:basic} gives $|\Rh(\psi;\bv,0)|\le\norm{\psi}_{\Wt}\big(\Ccal h^{-1}(h|\bk|)^{p+1}+\mathcal V_{\bk}\big)$.

\emph{Mass residual.} By \eqref{eq:orth}, $(\nabla\cdot\bu_h,q_h)_\Omega=0$ for all $q_h\in\Qh$, and trivially for constants, hence for all $q_h\in V_h^{p-1}$. Therefore, for $\psi\in L^2(I)$,
\[
|\Rh(\psi;\bm0,\phi_{\bk})| = \Big|\int_I\psi(\nabla\cdot\bu_h,\phi_{\bk}-\pi_Q\phi_{\bk})_\Omega\dd t\Big|\le\norm{\psi}_{L^2(I)}T^{1/2}\sup_{t}\norm{\nabla\cdot\bu_h}\,\norm{\phi_{\bk}-\pi_Q\phi_{\bk}},
\]
and the right-hand side is bounded by $C\norm{\psi}_{L^2(I)}T^{1/2}\E\,h^{-1}\Ms_{p}(h|\bk|)$ by the inverse inequality, \eqref{eq:energy} and \cref{lem:modes} with $p_Q=p-1$.
\end{proof}

\begin{remark}
The bound for $h|\bk|\ge1$, the ``trivial'' bound of the proof, equals $h^{-1}\le|\bk|$ and is attained, up to constants, when $h|\bk|\approx1$ (see \cref{sec:numerics}). In the envelope form of \eqref{eq:type} it reads $\rho_{\bk}\lesssim|\bk|$, which is the bound obtained in \cref{sec:model}. For the resolved scales, the momentum residual converges at the rate $h^{p}$, but the mass residual only at the rate $h^{p-1}$: the upper bound is hampered by the poor approximation properties of the pressure space and the lack of control of $\nabla\cdot\bu_h$. With the definition \eqref{eq:theta}, the resolution exponent of the method is $\theta=(p-1)/p$, i.e.\ $\theta=1/2$ for the lowest-order Taylor--Hood element.
\end{remark}

\subsection{Taylor--Hood elements with grad--div stabilisation}
A natural remedy for the limited resolved order of the mass residual is to choose the stabilising term $S_h$ in a way that gives additional control of $\nabla\cdot\bu_h$. Following \cite{franca1988two,olshanskii2004graddiv} we set
\[
S_h((\bu_h,p_h),(\bv_h,q_h)) := \gamma\,(\nabla\cdot\bu_h,\nabla\cdot\bv_h)_\Omega,\qquad\gamma>0 .
\]
The resulting scheme was advocated as a method suitable for high Reynolds number flow computations in \cite{john2010numerical}, and an error analysis with constants independent of the viscosity was proposed for smooth solutions in \cite{defrutos2018analysis}. By \eqref{eq:energy}, $\gamma\int_I\norm{\nabla\cdot\bu_h}^2\dd t\le\E^2$.

\begin{theorem}\label{thm:graddiv}
Let \cref{ass:fortin} hold. The residual spectrum of the Taylor--Hood method with grad--div stabilisation satisfies, for all $\bk\in\Z^d$, $h\in(0,1]$ and $\mu\ge0$,
\begin{align}
\rho^{\rm mom}_{\bk} &\le \Ccal\,(1+(\gamma T)^{1/2})\,h^{-1}\Ms_{p+1}(h|\bk|) + \mathcal V_{\bk}, \label{eq:gd_mom}\\
\rho^{\rm div}_{\bk} &\le C\,\gamma^{-1/2}\,\E\,\Ms_{p}(h|\bk|). \label{eq:gd_div}
\end{align}
In particular the method has residual spectrum of type $(1,p+1)$, with resolved order $\sigma-\kappa=p$, and the mass residual alone is of type $(0,p)$.
\end{theorem}
\begin{proof}
For $h|\bk|\ge1$ the proof of \cref{thm:th} applies verbatim. For $h|\bk|<1$, \eqref{eq:orth} now gives $\langle r_h,(\piF\bv,0)\rangle=\gamma(\nabla\cdot\bu_h,\nabla\cdot\piF\bv)_\Omega=\gamma(\nabla\cdot\bu_h,\nabla\cdot(\piF\bv-\bv))_\Omega$, where we used $\nabla\cdot\bv=0$; this is why only solenoidal Fourier modes are considered. Hence
\[
\Rh(\psi;\bv,0)=\int_I\psi\langle r_h,(\bv-\piF\bv,0)\rangle\dd t+\gamma\int_I\psi(\nabla\cdot\bu_h,\nabla\cdot(\piF\bv-\bv))_\Omega\dd t .
\]
The first term is bounded as in the proof of \cref{thm:th}. For the second,
\begin{multline*}
\Big|\gamma\int_I\psi(\nabla\cdot\bu_h,\nabla\cdot(\piF\bv-\bv))_\Omega\dd t\Big|\le\norm{\psi}_{L^\infty}(\gamma T)^{1/2}\Big(\gamma\int_I\norm{\nabla\cdot\bu_h}^2\dd t\Big)^{1/2}\sqrt d\,\norm{\nabla(\bv-\piF\bv)}\\ \le C\norm{\psi}_{L^\infty}(\gamma T)^{1/2}\E\,h^{-1}(h|\bk|)^{p+1}
\end{multline*}
by \eqref{eq:fortin_modes}. For the mass residual, the Galerkin orthogonality $(\nabla\cdot\bu_h,q_h)_\Omega=0$ for $q_h\in V_h^{p-1}$ still holds, and now
\[
|\Rh(\psi;\bm0,\phi_{\bk})|\le\norm{\psi}_{L^2(I)}\Big(\int_I\norm{\nabla\cdot\bu_h}^2\dd t\Big)^{1/2}\norm{\phi_{\bk}-\pi_Q\phi_{\bk}}\le C\norm{\psi}_{L^2(I)}\gamma^{-1/2}\E\,\Ms_p(h|\bk|).\qedhere
\]
\end{proof}

Compared to the standard Galerkin method the resolved order has improved by one, and the resolution exponent is $\theta=p/(p+1)$. Note also the trade-off in the choice of $\gamma$: the momentum bound grows like $\gamma^{1/2}$, whereas the bounds of the mass residual and of the divergence consistency term ($|\epsilon_h|\lesssim\gamma^{-1/2}$ by \cref{lem:Kh}) decrease like $\gamma^{-1/2}$. The behaviour of the unresolved scales is not improved: grad--div stabilisation does not act on the convective term, and the momentum residual still has the plateau $h^{-1}$ for $h|\bk|\ge1$.

\begin{remark}[An a posteriori improvement]\label{rem:gd_curl}
The plateau $h^{-1}$ in \cref{thm:th,thm:graddiv} comes only from the inverse inequality in the bound $|c(\bu_h;\bu_h,\bw)|\le C\norm{\bu_h}\norm{\nabla\bu_h}\norm{\bw}_{L^\infty}$. On the torus $\norm{\nabla\bv}^2=\norm{\nabla\times\bv}^2+\norm{\nabla\cdot\bv}^2$ for $\bv\in H^1(\Omega)^d$, so that for grad--div stabilisation
\[
\int_I\norm{\nabla\bu_h}\dd t\le\int_I\norm{\nabla\times\bu_h}\dd t+(T/\gamma)^{1/2}\E .
\]
If the computed vorticity satisfies $\int_I\norm{\nabla\times\bu_h}\dd t\le C_\omega h^{-1/2}$, the plateau and the first term in the proof of \cref{thm:graddiv} for $h|\bk|<1$ are therefore of order $h^{-1/2}$ (of order one if $\norm{\nabla\times\bu_h}$ is bounded). The second term, the consistency error of the grad--div term, is bounded by $C\gamma T^{1/2}\norm{\nabla\cdot\bu_h}_{Q}h^{-1}(h|\bk|)^{p+1}$; if in addition $\norm{\nabla\cdot\bu_h}_Q\le C_dh^{1/2}$, the momentum residual is of type $(\frac12,p+1)$, as for the CIP method of \cref{sec:sym}, and the bounds for the mass residual and for $\epsilon_h$ gain the factor $h^{1/2}$. Both conditions concern computed quantities and can be checked a posteriori. In two dimensions $\norm{\nabla\times\bu}$ is conserved by the Euler equations, but in three dimensions the vorticity is not controlled by the Euler equations, and for a field with a Kolmogorov spectrum resolved down to the scale $h$ one expects $\norm{\nabla\times\bu_h}\sim h^{-2/3}$, which violates the first condition.
\end{remark}

\section{Stabilised equal-order methods}\label{sec:stabilised}

In this section $p_Q=p\ge1$, so that the inf-sup condition fails and a pressure stabilisation is needed. Throughout, $U>0$ denotes a reference velocity, used to make the stabilisation parameters dimensionally consistent. We use the $L^2$-projection $\pih$ instead of a Fortin operator; since $\partial_t\bu_h(t)\in\Vh$ is $L^2$-orthogonal to $\bv-\pih\bv$, the time derivative drops out of the residual coefficients, and the pressure term is handled through the stabilisation.

\subsection{A nonconsistent streamline-diffusion method}\label{sec:supg}
We consider the stabilising term
\begin{equation}\label{eq:supg}
S_h(\bw;(\bu_h,p_h),(\bv_h,q_h)) := \big(\tau((\bw\cdot\nabla)\bu_h+\nabla p_h),(\bw\cdot\nabla)\bv_h+\nabla q_h\big)_\Omega + (\tau_d\nabla\cdot\bu_h,\nabla\cdot\bv_h)_\Omega,
\end{equation}
with $\tau=c_\tau h/U$ and $\tau_d=c_dhU$ for dimensionless constants $c_\tau,c_d>0$, and $\bw=\bu_h$ in \eqref{eq:fem}; in tables and figures we call it nonconsistent SUPG (streamline-upwind Petrov--Galerkin). This type of stabilisation was introduced for turbulent flows in \cite{hoffman2006new} as a nonconsistent variant of the velocity--pressure streamline-diffusion method of \cite{hansbo1990velocity}, claimed to be more suitable for the computation of turbulent flows. It is nonconsistent since the time derivative (and the viscous term) is not included in the residual in the first argument, so that $S_h$ does not vanish for the exact solution. In this case
\begin{equation}\label{eq:supg_norm}
\norm{(\bu_h\cdot\nabla)\bu_h+\nabla p_h}\le\tau^{-1/2}\snorm{(\bu_h,p_h)}_S,\qquad \norm{\nabla\cdot\bu_h}\le\tau_d^{-1/2}\snorm{(\bu_h,p_h)}_S .
\end{equation}
(If $\tau$ is limited by $c_\mu h^2/\mu$ for $\mu>Uh$, with a constant $c_\mu>0$, the resulting factor $\mu^{1/2}h^{-1}$ in the first bound gives a term of the form $\mathcal V_{\bk}$, as in the proof of \cref{thm:sym}, with a constant depending on $c_\mu$.)

\begin{theorem}\label{thm:supg}
The residual spectrum of the method \eqref{eq:fem}, \eqref{eq:supg} satisfies, for all $\bk\in\Z^d$, $h\in(0,1]$ and $\mu\ge0$,
\[
\rho^{\rm mom}_{\bk} \le \Ccal\,h^{-1/2}\Ms_{1}(h|\bk|) + \mathcal V_{\bk},\qquad
\rho^{\rm div}_{\bk} \le \Ccal\,h^{-1/2}\Ms_{1}(h|\bk|),
\]
where $\Ccal$ depends also on $U$, $c_\tau$ and $c_d$. In particular the method has residual spectrum of type $(\tfrac12,1)$, with resolved order $\sigma-\kappa=\frac12$ independently of $p$.
\end{theorem}
\begin{proof}
Let $\ba\in\mathcal A_{\bk}$, $\bv=\ba\phi_{\bk}$, $\bw:=\bv-\pih\bv$ and $\psi\in\Wt$. By \eqref{eq:orth},
\[
\langle r_h,(\bv,0)\rangle = \langle r_h,(\bw,0)\rangle + S_h(\bu_h;(\bu_h,p_h),(\pih\bv,0)).
\]
Since $(\partial_t\bu_h,\bw)_\Omega=0$ and $(p_h,\nabla\cdot\bw)_\Omega=-(\nabla p_h,\bw)_\Omega$,
\begin{equation}\label{eq:rw_stab}
\langle r_h,(\bw,0)\rangle = (\bff,\bw)_\Omega - ((\bu_h\cdot\nabla)\bu_h+\nabla p_h,\bw)_\Omega - \tfrac12((\nabla\cdot\bu_h)\bu_h,\bw)_\Omega - \mu(\nabla\bu_h,\nabla\bw)_\Omega .
\end{equation}
By \eqref{eq:supg_norm} and \cref{lem:modes},
\begin{multline*}
|((\bu_h\cdot\nabla)\bu_h+\nabla p_h,\bw)_\Omega|+\tfrac12|((\nabla\cdot\bu_h)\bu_h,\bw)_\Omega|\le\big(\tau^{-1/2}\norm{\bw}+\tfrac12\tau_d^{-1/2}\norm{\bu_h}\norm{\bw}_{L^\infty}\big)\snorm{(\bu_h,p_h)}_S\\
\le C(U^{1/2}+U^{-1/2}\E)\,h^{-1/2}\Ms_{p+1}(h|\bk|)\,\snorm{(\bu_h,p_h)}_S .
\end{multline*}
For the stabilisation term,
\begin{align*}
|S_h(\bu_h;(\bu_h,p_h),(\pih\bv,0))|&\le\snorm{(\bu_h,p_h)}_S\big(\tau^{1/2}\norm{\bu_h}\norm{\nabla\pih\bv}_{L^\infty}+\tau_d^{1/2}\norm{\nabla\cdot\pih\bv}\big)\\
&\le C(U^{-1/2}\E+U^{1/2})\,h^{-1/2}\Ms_1(h|\bk|)\,\snorm{(\bu_h,p_h)}_S,
\end{align*}
by \eqref{eq:m3}. The remaining terms in \eqref{eq:rw_stab} are bounded as in the proof of \cref{thm:th}. Multiplying by $\psi$, integrating in time and using $\int_I\snorm{(\bu_h,p_h)}_S\dd t\le T^{1/2}\E$ (by \eqref{eq:energy}) together with $\Ms_{p+1}\le\Ms_1$, we obtain the bound for $\rho^{\rm mom}_{\bk}$.

For the mass residual, \eqref{eq:orth} with $(\bv_h,q_h)=(\bm0,\pih\phi_{\bk})$ gives
\[
(\nabla\cdot\bu_h,\pih\phi_{\bk})_\Omega=-S_h(\bu_h;(\bu_h,p_h),(\bm0,\pih\phi_{\bk}))=-(\tau((\bu_h\cdot\nabla)\bu_h+\nabla p_h),\nabla\pih\phi_{\bk})_\Omega .
\]
Hence
\begin{align*}
|(\nabla\cdot\bu_h,\phi_{\bk})_\Omega|&\le\norm{\nabla\cdot\bu_h}\norm{\phi_{\bk}-\pih\phi_{\bk}}+\tau^{1/2}\snorm{(\bu_h,p_h)}_S\norm{\nabla\pih\phi_{\bk}}\\
&\le C\,U^{-1/2}h^{-1/2}\Ms_1(h|\bk|)\,\snorm{(\bu_h,p_h)}_S,
\end{align*}
and we conclude using $\int_I|\psi|\snorm{(\bu_h,p_h)}_S\dd t\le\norm{\psi}_{L^2(I)}\E$.
\end{proof}

We see that the resolution of the resolved scales is only $O(h^{1/2})$ regardless of the polynomial order. This is due to the lack of consistency of the stabilising term: the discrete approximation $\pih\bv$ of a smooth function is not small in the stabilisation semi-norm, since $\norm{\tau^{1/2}(\bu_h\cdot\nabla)\pih\bv}$ is in general of the order $h^{1/2}|\bk|\norm{\bu_h}$ for $h|\bk|<1$. The unresolved scales, on the other hand, are much better controlled than for the Galerkin methods: the plateau is $h^{-1/2}$ instead of $h^{-1}$ (envelope $|\bk|^{1/2}$ instead of $|\bk|$). The resolution exponent is $\theta=1/2$.

\subsection{Symmetric, weakly consistent stabilisation}\label{sec:sym}
Finally we discuss weakly consistent stabilised methods from the class of symmetric stabilisation methods. This class includes orthogonal subscale stabilisation \cite{codina2000stabilization,codina2002stabilized}, local projection stabilisation \cite{becker2004twolevel,braack2006local} and continuous interior penalty stabilisation \cite{douglas1976interior,burman2004edge}. The idea of all these methods is to use a scale separation into a space of coarse (resolved) scales and a space of fine (unresolved) scales, and to dissipate the component of the residual $(\bu_h\cdot\nabla)\bu_h+\nabla p_h$ that lies in the unresolved space, so that the dissipation is localised in frequency space and acts only on the highest frequencies. They have been advocated for the approximation of high Reynolds number and turbulent flows in several publications \cite{burman2007interior,burman2007continuous,guasch2013statistical,moura2022gradient}; however, the only theoretical justification so far has been an improved convergence of $O(h^{1/2})$, independent of the viscosity, under the condition that the solution is smooth.

We formulate the analysis under the following abstract assumptions and verify them for the continuous interior penalty method with $p=1$ below. We assume that the stabilisation splits as $S_h(\bw;(\bu_h,p_h),(\bv_h,q_h))=S_u(\bw;\bu_h,\bv_h)+S_p(p_h,q_h)$, with $S_u(\bw;\cdot,\cdot)$ and $S_p$ symmetric positive semi-definite, $S_p(q,1)=0$ for all $q$, and write $\snorm{\bv_h}_{S_u}$, $\snorm{q_h}_{S_p}$ for the associated semi-norms (with $\bw=\bu_h$). The term $\mu/h$ in (S1) allows stabilisation parameters that are reduced in the viscous regime $\mu>Uh$; for $\mu\le Uh$ it can be dropped.

\begin{assumption}\label{ass:sym}
There is a constant $C_S$, independent of $h$, $\mu$, $t$ and $\bk$, such that for all $t\in I$:
\begin{enumerate}[label=(S\arabic*),leftmargin=*]
\item (fluctuation control)
\begin{gather*}
\inf_{\bz_h\in\Vh}\norm{(\bu_h\cdot\nabla)\bu_h+\nabla p_h-\bz_h}\le C_S(U+\mu/h)^{1/2}h^{-1/2}\snorm{(\bu_h,p_h)}_S,\\ \inf_{y_h\in V_h^p}\norm{\nabla\cdot\bu_h-y_h}\le C_SU^{-1/2}h^{-1/2}\snorm{(\bu_h,p_h)}_S;
\end{gather*}
\item (weak consistency) for all $\bk\in\Z^d$ and $\ba\in\mathcal A_{\bk}$,
\begin{gather*}
\snorm{\pih(\ba\phi_{\bk})}_{S_u}\le C_S(U^{1/2}+U^{-1/2}\norm{\bu_h})\,h^{-1/2}\Ms_{p+1}(h|\bk|),\\ \snorm{\pih\phi_{\bk}}_{S_p}\le C_SU^{-1/2}h^{-1/2}\Ms_{p+1}(h|\bk|);
\end{gather*}
\item (inverse estimate) $\snorm{q_h}_{S_p}\le C_SU^{-1/2}h^{-1/2}\norm{q_h}$ for all $q_h\in V_h^p$.
\end{enumerate}
\end{assumption}

\begin{lemma}\label{lem:div}
Under \cref{ass:sym}, $\norm{\nabla\cdot\bu_h(t)}\le2C_SU^{-1/2}h^{-1/2}\snorm{(\bu_h,p_h)(t)}_S$.
\end{lemma}
\begin{proof}
Write $\nabla\cdot\bu_h=(\mathrm{Id}-\pih)\nabla\cdot\bu_h+\pih\nabla\cdot\bu_h$. The first term is bounded by (S1). For the second, the mass equation in \eqref{eq:fem} gives $(\nabla\cdot\bu_h,q_h)_\Omega=-S_p(p_h,q_h)$ for all $q_h\in\Qh$, and also for constants since both sides vanish; hence for all $q_h\in V_h^p$, by (S3),
$(\nabla\cdot\bu_h,q_h)_\Omega\le\snorm{p_h}_{S_p}\snorm{q_h}_{S_p}\le C_SU^{-1/2}h^{-1/2}\snorm{(\bu_h,p_h)}_S\norm{q_h}$, and $\norm{\pih\nabla\cdot\bu_h}=\sup_{q_h\in V_h^p}(\nabla\cdot\bu_h,q_h)_\Omega/\norm{q_h}$.
\end{proof}

\begin{theorem}\label{thm:sym}
Under \cref{ass:sym} the residual spectrum of the method satisfies, for all $\bk\in\Z^d$, $h\in(0,1]$ and $\mu\ge0$,
\[
\rho^{\rm mom}_{\bk}\le\Ccal\,h^{-1/2}\Ms_{p+1}(h|\bk|)+\mathcal V_{\bk},\qquad
\rho^{\rm div}_{\bk}\le\Ccal\,h^{-1/2}\Ms_{p+1}(h|\bk|),
\]
with $\Ccal$ depending also on $U$ and $C_S$. In particular the method has residual spectrum of type $(\tfrac12,p+1)$, with resolved order $\sigma-\kappa=p+\frac12$.
\end{theorem}
\begin{proof}
We proceed as in the proof of \cref{thm:supg}, starting from \eqref{eq:rw_stab}. Since $\bw=\bv-\pih\bv$ is $L^2$-orthogonal to $\Vh$, for every $\bz_h\in\Vh$
\begin{align*}
|((\bu_h\cdot\nabla)\bu_h+\nabla p_h,\bw)_\Omega| &= |((\bu_h\cdot\nabla)\bu_h+\nabla p_h-\bz_h,\bw)_\Omega|\\
&\le C_S(U^{1/2}h^{-1/2}+\mu^{1/2}h^{-1})\snorm{(\bu_h,p_h)}_S\norm{\bw}
\end{align*}
by (S1). By \cref{lem:div}, $\tfrac12|((\nabla\cdot\bu_h)\bu_h,\bw)_\Omega|\le C_SU^{-1/2}h^{-1/2}\snorm{(\bu_h,p_h)}_S\norm{\bu_h}\norm{\bw}_{L^\infty}$, and by (S2)
\[
|S_u(\bu_h;\bu_h,\pih\bv)|\le\snorm{\bu_h}_{S_u}\snorm{\pih\bv}_{S_u}\le C_S(U^{1/2}+U^{-1/2}\norm{\bu_h})\,h^{-1/2}\Ms_{p+1}(h|\bk|)\,\snorm{(\bu_h,p_h)}_S .
\]
Using \eqref{eq:m1}, integrating in time against $\psi$ and using $\int_I\snorm{(\bu_h,p_h)}_S\dd t\le T^{1/2}\E$ yields the bound for $\rho^{\rm mom}_{\bk}$; the contribution of $\mu^{1/2}h^{-1}$ is of the form $\mathcal V_{\bk}$ (with a constant depending on $C_S$), since $h^{-1}\Ms_{p+1}(h|\bk|)\le|\bk|\Ms_p(h|\bk|)$. For the mass residual, the mass equation gives $(\nabla\cdot\bu_h,\pih\phi_{\bk})_\Omega=-S_p(p_h,\pih\phi_{\bk})$, so that
\[
|(\nabla\cdot\bu_h,\phi_{\bk})_\Omega|\le\norm{\nabla\cdot\bu_h}\norm{\phi_{\bk}-\pih\phi_{\bk}}+\snorm{p_h}_{S_p}\snorm{\pih\phi_{\bk}}_{S_p}\le C\,U^{-1/2}h^{-1/2}\Ms_{p+1}(h|\bk|)\,\snorm{(\bu_h,p_h)}_S,
\]
by \cref{lem:div}, \eqref{eq:m1} and (S2), and we conclude as before.
\end{proof}

\paragraph{The continuous interior penalty method.}
We now consider the continuous interior penalty (CIP) method \cite{douglas1976interior}, analysed for the equations of fluid mechanics in \cite{burman2006continuous,burman2007continuous} and applied to under-resolved turbulent flow in \cite{moura2022gradient}. Dissipation is introduced through a penalty on the jump of the gradient over interior faces, and the divergence is controlled by a bulk grad--div term as in \eqref{eq:supg}:
\begin{equation}\label{eq:cip}
\begin{aligned}
S_u(\bw;\bu_h,\bv_h)&:=\gamma_u\sum_{F\in\mathcal F_h}h_F^2\int_F\omega_F(\bw)\,\jump{\partial_{\bn}\bu_h}\cdot\jump{\partial_{\bn}\bv_h}\dd s+(\tau_d\nabla\cdot\bu_h,\nabla\cdot\bv_h)_\Omega,\\
S_p(p_h,q_h)&:=\gamma_p\sum_{F\in\mathcal F_h}h_F^2\min\{U^{-1},h_F/\mu\}\int_F\jump{\partial_{\bn}p_h}\,\jump{\partial_{\bn}q_h}\dd s,
\end{aligned}
\end{equation}
with parameters $\gamma_u,\gamma_p>0$, $\bw=\bu_h$, $\tau_d=c_dhU$, where $h_F$ is the diameter of the face $F$, $\bn$ a unit normal to $F$, $\partial_{\bn}$ the normal derivative, and
\begin{equation}\label{eq:cip_weight}
\omega_F(\bw):=U^{-1}|\bw\cdot\bn|^2+\varepsilon\,U,\qquad 0\le\varepsilon\le1 .
\end{equation}
Taking $\varepsilon>0$ adds some crosswind diffusion, which is known to enhance control of the divergence of the solution, allowing for $\tau_d=0$, but smears sharp layers. The pressure weight is that of \cite{burman2007continuous}; it equals $U^{-1}$ if $\mu\le Uh_F$, in particular for $\mu=0$. Since $\bu_h$ and $p_h$ are continuous, only their normal derivatives jump: $\jump{\nabla\bu_h}=\jump{\partial_{\bn}\bu_h}\otimes\bn$ and $\jump{\nabla p_h}=\jump{\partial_{\bn}p_h}\bn$.

\begin{lemma}\label{lem:cip}
For $p=1$ the CIP method \eqref{eq:cip} satisfies \cref{ass:sym}, with $C_S$ depending only on $d$, the mesh regularity, $\gamma_u$, $\gamma_p$ and $c_d$, for every $\varepsilon\in[0,1]$.
\end{lemma}
\begin{proof}
We use the Oswald quasi-interpolation estimate \cite{karakashian2003posteriori,burman2007continuoushp}: for every piecewise polynomial $g$ of degree at most $p$ there is $I_{\rm os}g\in V_h^p$ with
\begin{equation}\label{eq:oswald}
\norm{g-I_{\rm os}g}^2\le C\sum_{F\in\mathcal F_h}h_F\norm{\jump{g}}^2_{L^2(F)} .
\end{equation}
(S1): For $p=1$, $\bm g:=(\bu_h\cdot\nabla)\bu_h+\nabla p_h$ is piecewise linear with $\jump{\bm g}=(\bu_h\cdot\bn)\jump{\partial_{\bn}\bu_h}+\jump{\partial_{\bn}p_h}\bn$. Since $|\bu_h\cdot\bn|^2\le U\omega_F(\bu_h)$ and $\max\{U,\mu/h_F\}\min\{U^{-1},h_F/\mu\}=1$, \eqref{eq:oswald} and $h\le Ch_F$ give
\begin{align*}
h\norm{\bm g-I_{\rm os}\bm g}^2&\le C\sum_{F}h_F^2\big(\norm{(\bu_h\cdot\bn)\jump{\partial_{\bn}\bu_h}}^2_{L^2(F)}+\norm{\jump{\partial_{\bn}p_h}}^2_{L^2(F)}\big)\\
&\le C(U+\mu/h)\big(\gamma_u^{-1}\snorm{\bu_h}_{S_u}^2+\gamma_p^{-1}\snorm{p_h}_{S_p}^2\big).
\end{align*}
The second inequality in (S1) holds even with $y_h=0$, since $\norm{\nabla\cdot\bu_h}\le\tau_d^{-1/2}\snorm{\bu_h}_{S_u}=(c_dhU)^{-1/2}\snorm{\bu_h}_{S_u}$.
(S2): Let $\bv=\ba\phi_{\bk}$. For the face term we use $\omega_F(\bu_h)\le U^{-1}|\bu_h|^2+U$, $\norm{\jump{\partial_{\bn}\pih\bv}}_{L^\infty(F)}\le\norm{\jump{\nabla\pih\bv}}_{L^\infty(F)}$, the discrete trace inequality $\sum_F\norm{\bu_h}^2_{L^2(F)}\le Ch^{-1}\norm{\bu_h}^2$ and $\sum_F|F|\le Ch^{-1}$, with $|F|$ the $(d-1)$-dimensional measure of $F$:
\[
\gamma_u\sum_Fh_F^2\int_F\omega_F(\bu_h)|\jump{\partial_{\bn}\pih\bv}|^2\dd s\le Ch\,(U^{-1}\norm{\bu_h}^2+U)\max_F\norm{\jump{\nabla\pih\bv}}^2_{L^\infty(F)},
\]
and \eqref{eq:m4} gives the bound of (S2). For the grad--div term, $\nabla\cdot\pih\bv=\nabla\cdot(\pih\bv-\bv)$ since $\bv$ is solenoidal, so that by \eqref{eq:m2} and \eqref{eq:m3} $\tau_d^{1/2}\norm{\nabla\cdot\pih\bv}\le C(hU)^{1/2}\min\{|\bk|\Ms_p(h|\bk|),h^{-1}\}\le CU^{1/2}h^{-1/2}\Ms_{p+1}(h|\bk|)$. The bound for $S_p$ is obtained in the same way as for the face term, since its weight is at most $U^{-1}$.
(S3): Since the weight is at most $U^{-1}$, the discrete trace and inverse inequalities give $\snorm{q_h}_{S_p}^2\le C\gamma_pU^{-1}h^2\sum_K h^{-1}\norm{\nabla q_h}^2_{L^2(K)}\le CU^{-1}h^{-1}\norm{q_h}^2$.
\end{proof}

Thus the face weight only needs to control the normal component of the convecting velocity; the part $\varepsilon U$ is used only in \eqref{eq:Kh_cip}. If $\varepsilon>0$, the second inequality in (S1) also holds with $\tau_d=0$: $\nabla\cdot\bu_h$ is piecewise polynomial with $\jump{\nabla\cdot\bu_h}=\jump{\partial_{\bn}\bu_h}\cdot\bn$, so that \eqref{eq:oswald} and the part $\varepsilon U$ of the weight give it with $C_S$ depending on $\varepsilon^{-1/2}$.

\paragraph{Higher polynomial degrees.}
For $p\ge2$ the proof of \cref{lem:cip} still gives (S2), (S3), and (S1) except for the convective term; hence \cref{lem:div}, the bound for $\rho^{\rm div}_{\bk}$ in \cref{thm:sym} and, for $\varepsilon>0$, \eqref{eq:Kh_cip} hold for all $p$. The convective part of (S1), however, does not follow from the stabilisation: as observed in \cite[Sect.~2]{moura2022gradient}, the gradient jumps control the jumps of $(\bu_h\cdot\nabla)\bu_h$, but not its part $(\mathrm{Id}-\pi^X)(\bu_h\cdot\nabla)\bu_h$ of degree $2p-1$, where $\pi^X$ is the elementwise $L^2$-projection onto $[\mathbb P_p(K)]^d$ (the superscript $X$ refers to the broken, elementwise polynomial space). For the method as it stands we have the following bound.

\begin{proposition}\label{prop:cip_highp}
Let $p\ge2$ and let $\Pi_1^X$ denote the elementwise $L^2$-projection onto $[\mathbb P_1(K)]^d$. The residual spectrum of the CIP method \eqref{eq:cip} satisfies, for all $\bk\in\Z^d$, $h\in(0,1]$ and $\mu\ge0$, the bound for $\rho^{\rm div}_{\bk}$ of \cref{thm:sym} and
\[
\rho^{\rm mom}_{\bk}\le\Ccal\,h^{-1/2}\Ms_{p+1}(h|\bk|)+C\,h^{-1}\Ms_{p+1}(h|\bk|)\int_I\norm{(\mathrm{Id}-\Pi_1^X)\bu_h}^2\dd t+\mathcal V_{\bk}.
\]
In particular it is of type $(1,p+1)$, and of type $(\frac12,p+1)$ if $\int_I\norm{(\mathrm{Id}-\Pi_1^X)\bu_h}^2\dd t\le Ch^{1/2}$.
\end{proposition}
\begin{proof}
Only the bound of $(\bm g,\bw)_\Omega$ in the proof of \cref{thm:sym}, $\bm g:=(\bu_h\cdot\nabla)\bu_h+\nabla p_h$, $\bw=\bv-\pih\bv$, changes. Let $\bm g_1:=(\mathrm{Id}-\pi^X)\bm g$, $\bm g_2:=\pi^X\bm g$ and $\bm z:=(\mathrm{Id}-\Pi_1^X)\bu_h$. Since $\nabla\Pi_1^X\bu_h$ is elementwise constant, $(\bu_h\cdot\nabla)\bu_h-(\bm z\cdot\nabla)\bm z$ and $\nabla p_h$ are elementwise in $[\mathbb P_p]^d$, so that $\bm g_1=(\mathrm{Id}-\pi^X)(\bm z\cdot\nabla)\bm z$. By norm equivalence on $\mathbb P_{2p-1}(K)$ and an inverse inequality, $\norm{\bm g_1}_{L^1(K)}\le C\norm{\bm z}_{L^2(K)}\norm{\nabla\bm z}_{L^2(K)}\le Ch^{-1}\norm{\bm z}_{L^2(K)}^2$. Since $I_{\rm os}\bm g_2\in\Vh$ is orthogonal to $\bw$, $(\bm g,\bw)_\Omega=(\bm g_1,\bw)_\Omega+(\bm g_2-I_{\rm os}\bm g_2,\bw)_\Omega$. The local argument behind \eqref{eq:oswald} also gives $\norm{\bm g_2-I_{\rm os}\bm g_2}_{L^1}\le C\sum_Fh_F\norm{\jump{\bm g_2}}_{L^1(F)}$, and $\jump{\bm g_2}=\jump{\bm g}-\jump{\bm g_1}$. By a discrete trace inequality, $h_F\norm{\bm g_1|_K}_{L^1(F)}\le C\norm{\bm g_1}_{L^1(K)}$, and by the Cauchy--Schwarz inequality, $\sum_F|F|h_F\le C$ and the proof of \cref{lem:cip}, $\sum_Fh_F\norm{\jump{\bm g}}_{L^1(F)}\le C(U+\mu/h)^{1/2}h^{-1/2}\snorm{(\bu_h,p_h)}_S$. Hence
\[
|(\bm g,\bw)_\Omega|\le C\norm{\bw}_{L^\infty}\big((U+\mu/h)^{1/2}h^{-1/2}\snorm{(\bu_h,p_h)}_S+h^{-1}\norm{(\mathrm{Id}-\Pi_1^X)\bu_h}^2\big),
\]
and we conclude with \eqref{eq:m1} as in the proof of \cref{thm:sym}. The type $(1,p+1)$ follows from $\norm{(\mathrm{Id}-\Pi_1^X)\bu_h}\le\norm{\bu_h}\le\E$.
\end{proof}

The integral in \cref{prop:cip_highp}, the energy of $\bu_h$ in its elementwise content of degree $2$ to $p$, is computable a posteriori; snapshot values are reported in \cref{sec:num2d}.

\begin{remark}[Elementwise stabilisation]\label{rem:cip_elementwise}
Adding to $S_u$, with a parameter $\gamma_K>0$ and $(\cdot,\cdot)_K$ the $L^2(K)$ inner product, the elementwise term $\gamma_K\sum_K hU^{-1}\big((\mathrm{Id}-\pi^X)(\bw\cdot\nabla)\bu_h,(\mathrm{Id}-\pi^X)(\bw\cdot\nabla)\bv_h\big)_K$ (a variant of that of \cite{moura2022gradient}), which vanishes for $p=1$, restores (S1) for all $p$, since $\inf_{\bz_h\in\Vh}\norm{\bm g-\bz_h}\le C\big(\norm{\bm g_1}+(\sum_Fh_F\norm{\jump{\bm g}}^2_{L^2(F)})^{1/2}\big)$ \cite[Prop.~1]{moura2022gradient}. For smooth solutions it is weakly consistent of optimal order, $\norm{h^{1/2}(\mathrm{Id}-\pi^X)(\bu\cdot\nabla)\bu}\le Ch^{p+1/2}|(\bu\cdot\nabla)\bu|_{H^p(\Omega)}$, as for local projection stabilisation. The weak consistency (S2) used here, which concerns the projected test function with the discrete convecting velocity and uses no regularity of $\bu_h$, is lower: as in the proof above, $(\mathrm{Id}-\pi^X)(\bu_h\cdot\nabla)\pih\bv=(\mathrm{Id}-\pi^X)(\nabla\pih\bv-\Pi_0^X\nabla\pih\bv)(\bu_h-\Pi_1^X\bu_h)$, $\Pi_0^X$ the elementwise mean, so that, with only $\norm{\bu_h-\Pi_1^X\bu_h}\le\norm{\bu_h}$, (S2) holds with $\Ms_2$ in place of $\Ms_{p+1}$ and the type is $(\frac12,2)$ for all $p$ and hence gives no improvement here compared to the results for $p=1$, without the term. For local projection and orthogonal subscale stabilisation (S1) is built into the method (when the fine-scale space contains the convective term) and the analysis hinges on (S2).
\end{remark}

\subsection{Negative-norm versions of the bounds}\label{sec:negative_methods}
The per-mode bounds above were obtained by testing with $L^\infty$-functions; this is what allowed us to use only the $L^2$-bound \eqref{eq:energy} on $\bu_h$. To obtain bounds of $\ell^2$-type (\cref{def:l2type}), the test functions are measured in $L^2$-based norms, and we need in addition that the discrete velocities are bounded,
\begin{equation}\label{eq:Uinf}
\norm{\bu_h}_{L^\infty(Q)}\le U_\infty ,
\end{equation}
with $U_\infty$ independent of $h$ and $\mu$. This is physically reasonable for turbulent flows (velocities are bounded, gradients are not) and can be checked a posteriori, but does not follow from \eqref{eq:energy}. For the Taylor--Hood element we also use the $L^2$-version of \cref{ass:fortin}: $\norm{\bv-\piF\bv}+h\norm{\nabla(\bv-\piF\bv)}\le C_Fh^s|\bv|_{H^s}$ for $1\le s\le p+1$, which holds for $p=2$ by \cite[Prop.~2.5]{diening2022fortin} and is assumed for $p\ge3$. Since the test functions are now only in $L^2$-based spaces, we extend the residual \eqref{eq:residual} (for $\mu=0$) to $\bv\in[L^2(\Omega)]^d$ by writing the pressure term as $-(\nabla p_h,\bv)_\Omega$, which is possible since $p_h\in H^1(\Omega)$, and the convective term as $-(\bu_h\cdot\nabla\bu_h+\frac12(\nabla\cdot\bu_h)\bu_h,\bv)_\Omega$. For symmetric stabilisations we replace (S2) by
\begin{enumerate}[label=(S2$'$),leftmargin=*]
\item $\snorm{\pih\bv}_{S_u}\le C_Sh^{s-1/2}\norm{\bv}_{H^s}$ for solenoidal $\bv$ and $\snorm{\pih q}_{S_p}\le C_Sh^{s-1/2}\norm{q}_{H^s}$, for $s\in\{0,p+1\}$.
\end{enumerate}
For the CIP method with $p=1$ and under \eqref{eq:Uinf}, (S2$'$) follows from the discrete trace and inverse inequalities for $s=0$, and for $s=2$ from $\jump{\partial_{\bn}\pih\bv}=\jump{\partial_{\bn}(\pih\bv-\bv)}$, the continuous trace inequality $\norm{z}_{L^2(F)}^2\le C(h^{-1}\norm{z}_{L^2(K)}^2+h\norm{\nabla z}_{L^2(K)}^2)$ applied to $z=\nabla(\pih\bv-\bv)$, and the approximation properties of $\pih$; for the bulk term in \eqref{eq:cip}, $\tau_d^{1/2}\norm{\nabla\cdot\pih\bv}$ is bounded by $Ch^{-1/2}\norm{\bv}$ by an inverse inequality, and by $Ch^{1/2}\norm{\nabla(\pih\bv-\bv)}\le Ch^{3/2}\norm{\bv}_{H^2}$ since $\bv$ is solenoidal; the same argument applies for $p\ge2$. For CIP with $p\ge2$, the proof of \cref{prop:cip_highp} with $\norm{\bm g_1}\le CU_\infty h^{-1}\norm{(\mathrm{Id}-\Pi_1^X)\bu_h}$ gives $\ell^2$-type $(1,p+1)$, and $(\frac12,p+1)$ if $\int_I\norm{(\mathrm{Id}-\Pi_1^X)\bu_h}\dd t\le Ch^{1/2}$.

\begin{proposition}\label{prop:negative}
Let $\mu=0$ and assume \eqref{eq:Uinf} and, for the Taylor--Hood element, the $L^2$-version of \cref{ass:fortin}. Then the residuals of the Taylor--Hood Galerkin method, the Taylor--Hood method with grad--div stabilisation, the nonconsistent streamline-diffusion method and the symmetric stabilisation satisfying (S1), (S2$'$), (S3) are of $\ell^2$-type $(1,p)$, $(1,p+1)$, $(\tfrac12,1)$ and $(\tfrac12,p+1)$, respectively; for the Taylor--Hood Galerkin method the momentum part alone is of $\ell^2$-type $(1,p+1)$. The constants depend on $U_\infty$ in addition to the quantities listed in \cref{sec:fem}.
\end{proposition}
\begin{proof}[Sketch of proof]
The proofs follow those of \cref{thm:th,thm:graddiv,thm:supg,thm:sym}, with $(\bv,q)$ now time-dependent; we only indicate the modifications. All $L^\infty$-norms of the test function are replaced by $L^2$-norms, the price being a factor $U_\infty$ from \eqref{eq:Uinf} in the convective terms, e.g.\ $|c(\bu_h;\bu_h,\bw)|\le CU_\infty h^{-1}\norm{\bu_h}\norm{\bw}$, and in the weights $|\bu_h|^2$ in the CIP stabilisation and $\bu_h$ in \eqref{eq:supg}. \Cref{lem:modes} is replaced by the standard estimates $\norm{\bv-\pih\bv}\le Ch^s|\bv|_{H^s}$ ($0\le s\le p+1$), $\norm{\nabla(\bv-\pih\bv)}\le Ch^{s-1}|\bv|_{H^s}$ ($1\le s\le p+1$) and $\norm{\nabla\pih\bv}\le C\min\{\norm{\nabla\bv},h^{-1}\norm{\bv}\}$, valid on quasi-uniform meshes \cite{crouzeix1987stability}, and by their analogues for $\piF$ and $\pi_Q$. In the Taylor--Hood case the time derivative is integrated by parts, with $\partial_t\bw=(\mathrm{Id}-\piF)\partial_t\bv$. For the first bound in \eqref{eq:l2type} ($s=0$) no interpolant is needed for the Galerkin methods: the bound follows from stability alone, as in the case $h|\bk|\ge1$ of \cref{thm:th}, and gives $Ch^{-1}\norm{(\bv,q)}_{X_0}$. For the stabilised methods we still split $\bv=\pih\bv+(\bv-\pih\bv)$ and use (S1) or \eqref{eq:supg_norm} and (S2$'$) with $s=0$, which gives $Ch^{-1/2}\norm{(\bv,q)}_{X_0}$. For the second bound, the same steps with $s=\sigma$ give $Ch^{\sigma-\kappa}\norm{(\bv,q)}_{X_\sigma}$; for instance, for the nonconsistent method the term $(\tau((\bu_h\cdot\nabla)\bu_h+\nabla p_h),(\bu_h\cdot\nabla)\pih\bv)_\Omega$ is bounded by $\tau^{1/2}U_\infty\norm{\nabla\pih\bv}\snorm{(\bu_h,p_h)}_S\le Ch^{1/2}\norm{\bv}_{H^1}\snorm{(\bu_h,p_h)}_S$, which limits the resolved order to $\sigma-\kappa=\frac12$, as in \cref{thm:supg}.
\end{proof}

\subsection{Summary}
The results of \cref{sec:galerkin,sec:stabilised} are summarised in \cref{tab:summary}. The best resolved order is obtained for the symmetric stabilisation, which shares the slowest growth of the unresolved coefficients with the nonconsistent method; its resolution exponent is below the value $\theta=1$ of a method with $\kappa=0$ only by $\frac1{2(p+1)}$.

\begin{table}[ht]
\centering
\caption{Residual spectrum types $(\kappa,\sigma)$ of \eqref{eq:type}: resolved order $\sigma-\kappa$, resolution exponent $\theta=(\sigma-\kappa)/\sigma$ and the adjoint regularity required for convergence, for per-mode adjoint decay $(1+|\bk|)^{-\beta}$ in \cref{thm:abstract} and, under \eqref{eq:Uinf}, for adjoints in $X_s$ in \cref{thm:negative}; these concern the residual part of the error, without $\epsilon_h$. For Taylor--Hood $p\ge2$; for CIP with $p\ge2$ the entries hold under the condition of \cref{prop:cip_highp} and its $\ell^2$-analogue in \cref{sec:negative_methods}, and $(\kappa,\sigma)=(1,p+1)$ in general. Inviscid case. Last row: the benchmark type $(0,p+1)$ of \cref{sec:burgers}, with residual coefficients of the size of the projection error \eqref{eq:approx1d} of a Fourier mode (a reference value: no lower bound for the residual is proved).}\label{tab:summary}
\begin{tabular}{lccccccc}
\toprule
Method & spaces & $\kappa$ & $\sigma$ & $\sigma-\kappa$ & $\theta$ & $\beta>$ & $s>$\\
\midrule
Galerkin (TH), momentum & $\mathbb P_p/\mathbb P_{p-1}$ & $1$ & $p+1$ & $p$ & $\frac{p}{p+1}$ & $d+1$ & $1$\\
Galerkin (TH), mass & & $1$ & $p$ & $p-1$ & $\frac{p-1}{p}$ & $d+1$ & $1$\\
TH + grad--div & $\mathbb P_p/\mathbb P_{p-1}$ & $1$ & $p+1$ & $p$ & $\frac{p}{p+1}$ & $d+1$ & $1$\\
nonconsistent SUPG & $\mathbb P_p/\mathbb P_p$ & $\frac12$ & $1$ & $\frac12$ & $\frac12$ & $d+\frac12$ & $\frac12$\\
symmetric (CIP) & $\mathbb P_p/\mathbb P_p$ & $\frac12$ & $p+1$ & $p+\frac12$ & $\frac{2p+1}{2p+2}$ & $d+\frac12$ & $\frac12$\\
\midrule
benchmark ($L^2$-projection) & & $0$ & $p+1$ & $p+1$ & $1$ & $d$ & $0$\\
\bottomrule
\end{tabular}
\end{table}

\section{Numerical examples}\label{sec:numerics}

We compute residual spectra for the model problems of \cref{sec:model} and for an under-resolved two-dimensional flow, and compare them to the errors in the low Fourier modes. The residual is evaluated for the fully discrete scheme, with the time derivative treated as in \eqref{eq:residual_tw}: for the weight $\psi\equiv1$ on a window $[t_0,t_1]$ (that is, with $I$ replaced by $(t_0,t_1)$ and the error at $t_0$ as initial error),
\begin{equation}\label{eq:discrete_residual}
\Rh(1;\bv,q) = -(\bu_h(t_1)-\bu_h(t_0),\bv)_\Omega + \sum_{n}\Delta t\,\big\langle \tilde r_h^{\,n},(\bv,q)\big\rangle,
\end{equation}
where the sum runs over the time steps of length $\Delta t$ in $[t_0,t_1]$ and $\tilde r_h^{\,n}$ is the spatial part of \eqref{eq:residual} evaluated with the quadrature in time of the time-stepping scheme. Then \eqref{eq:orth} holds exactly at the fully discrete level. We expect the bounds of \cref{sec:galerkin,sec:stabilised} to carry over to this algebraic residual, with the convecting velocity $\bw^n$ (an extrapolated velocity that will be defined for each case below) of the scheme in place of $\bu_h$ in the convective terms and in the stabilisation. \Cref{prop:representation}, which holds for every $\bu_h\in W^{1,1}(I;\Vh)$, applies instead to the piecewise linear interpolant in time of the discrete solution, with residual \eqref{eq:residual}.
The Fourier coefficients of the piecewise polynomial fields in \eqref{eq:discrete_residual} are computed to rounding accuracy, well beyond the mesh resolution, by Gauss quadrature on each element combined with fast Fourier transforms.

\subsection{One space dimension}\label{sec:num1d}

\begin{figure}[!htb]
\centering
\includegraphics[width=\textwidth]{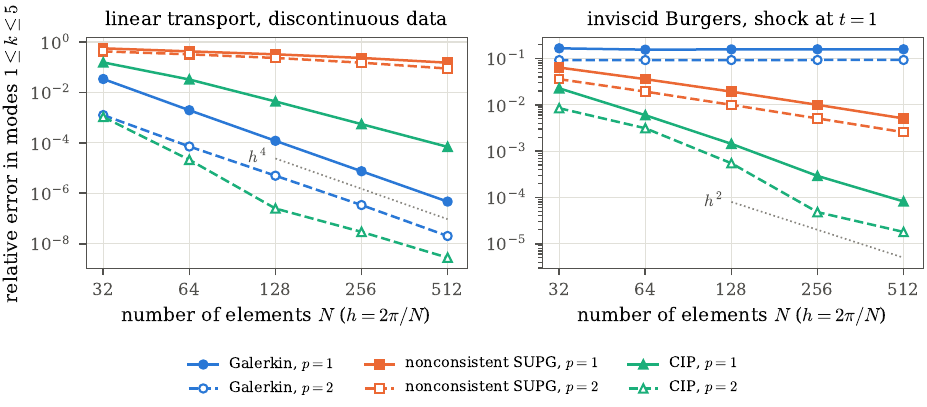}
\caption{Relative error in the Fourier modes $1\le k\le5$ at the final time for the transport equation with discontinuous data (left) and the inviscid Burgers equation after shock formation (right). Filled markers: $p=1$; open markers: $p=2$; dotted: reference slopes $h^4$ (left) and $h^2$ (right).}\label{fig:burgers_errors}
\end{figure}

\begin{figure}[!htb]
\centering
\includegraphics[width=\textwidth]{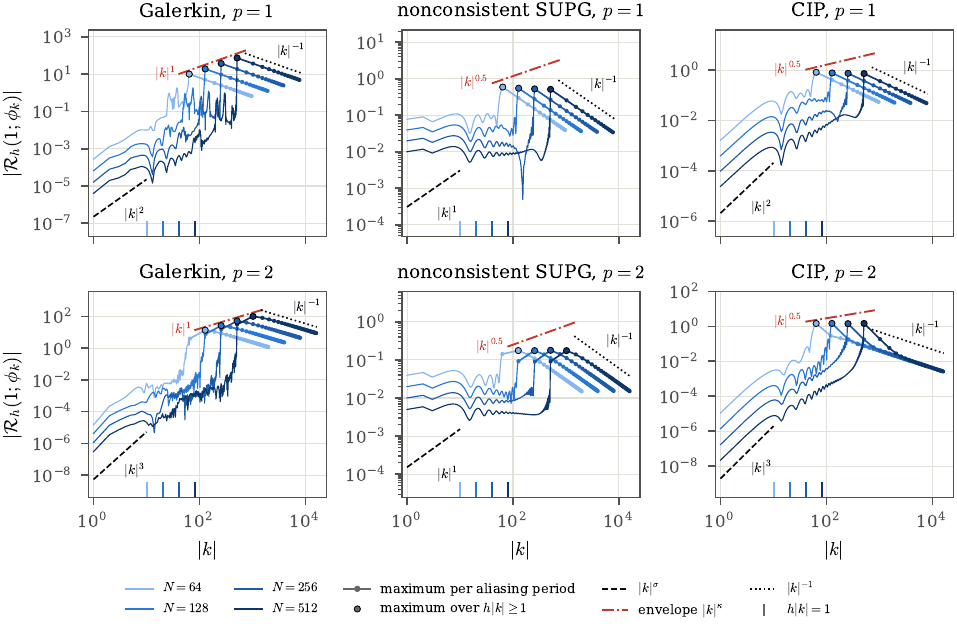}
\caption{Time-weak residual spectra $|\Rh(1;\phi_k)|$ ($\psi\equiv1$ on $[0,2]$) for the inviscid Burgers equation against $|k|\le16pN$, for $N=64,\dots,512$ (light to dark): the full spectrum up to the first aliasing peak, beyond it the maxima over the aliasing periods $|k-mN|<\frac N2$ (dots); circles: maximum over $h|k|\ge1$; ticks on the axis: the resolution limit $h|k|=1$ for each $N$. Dashed: growth $|k|^\sigma$ of the bounds ($\sigma=p+1$ for the Galerkin and CIP methods, $\sigma=1$ for SUPG); dash-dotted: their envelope $|k|^\kappa$; dotted: $|k|^{-1}$. The CIP method has the optimal spectrum, of type $(0,p+1)$.}\label{fig:burgers_spectra}
\end{figure}

We use uniform meshes of $\T$ with $N\in\{32,64,\dots,512\}$ elements, $h=2\pi/N$, and continuous elements of degree $p\in\{1,2\}$. Besides the Galerkin methods \eqref{eq:gal1d} and \eqref{eq:galburgers} we consider the one-dimensional analogues of the methods of \cref{sec:stabilised}: the nonconsistent streamline-diffusion method, $S_h(u_h;u_h,v_h)=(\tau w\partial_xu_h,w\partial_xv_h)_\Omega$ with $\tau=h/(2pU)$ and the convecting velocity $w=1$ (transport) or $w=u_h$ (Burgers), and the CIP method, $S_h(u_h;u_h,v_h)=\gamma_u\sum_{i}h^2\,\omega_i\,\jump{\partial_xu_h}(x_i)\jump{\partial_xv_h}(x_i)$ over the mesh nodes, with $\omega_i=U$ (transport) or $\omega_i=U^{-1}|u_h(x_i)|^2$ (Burgers), $\gamma_u=0.05$ for $p=1$ and $\gamma_u=0.01$ for $p=2$. With $r_h$ of \cref{sec:model} and as in \eqref{eq:residual_tw} and \cref{def:spectrum}, $\Rh(\psi;v):=\int_I\psi\langle r_h(t),v\rangle\dd t$ and $\rho_k:=\sup_\psi|\Rh(\psi;\phi_k)|/\norm{\psi}_{\Wt}$; we compute $|\Rh(1;\phi_k)|$. For the Burgers equation the proofs of \cref{thm:supg,thm:sym} apply with obvious modifications and give types $(\frac12,1)$ and $(\frac12,p+1)$ (for $p=2$ under the condition of \cref{prop:cip_highp}), while \eqref{eq:rkburgers} gives type $(1,p+1)$ for the Galerkin method. We integrate in time with the classical fourth-order Runge--Kutta method, with $\Delta t\le0.2h/(pU)$ and small enough for the stabilised operators to be inside the stability region.
Dividing $\Delta t$ by about four changes the residual coefficients for $1\le k\le8$ and the maxima in \cref{fig:burgers_spectra} by less than $0.4\%$ for the Galerkin method and $10^{-4}$ for the stabilised methods, and the errors in \cref{fig:burgers_errors} by less than $0.4\%$, except for the transport errors below $10^{-6}$ of the CIP method with $p=2$ (up to $6\%$); the time discretisation therefore plays no visible role in this section. For the transport equation $u_0(x)=\operatorname{sign}(\sin x)$, $U=1$, $T=2\pi$; for the Burgers equation $\mu=0$, $f=0$, $u_0(x)=\frac12+\sin x$, $U=\frac32$ and $T=2$, so that a shock forms at $t=1$. The exact entropy solution is obtained by a Galilean transformation from that for $u_0=\sin x$.

\Cref{fig:burgers_errors} shows the relative error in the Fourier modes $1\le k\le5$ at the final time, $(\sum_{k=1}^5|\hat e_k(T)|^2/\sum_{k=1}^5|\hat u_k(T)|^2)^{1/2}$. For the linear problem all methods converge in the low modes in spite of the discontinuous data, as predicted by \eqref{eq:ek1d} (the Galerkin method even superconverges on these uniform meshes). For the Burgers equation the Galerkin method, which conserves energy, does not converge in the low modes, the nonconsistent method converges at first order and the CIP method at second to third order. The residual spectra in \cref{fig:burgers_spectra} and the exponents in \cref{tab:rates1d} show why; they test the bounds in $h$ and in $k$ separately. \emph{Dependence on $h$.} The spectra show the three regimes of \cref{sec:burgers}. The peaks, at the first aliasing wave number $k=N$ (at $k=2N$ for the Galerkin and SUPG methods with $p=2$), grow like $h^{-0.96}$ for the Galerkin method, as $h^{-\kappa}=h^{-1}$ in \eqref{eq:rkburgers}, and lie on the envelope $|k|$, whereas those of the stabilised methods are bounded uniformly in $h$ (theoretical upper bound $h^{-1/2}$). Beyond the peaks the maxima over the aliasing periods $|k-mN|<\frac N2$ decay like $|k|^{-1}$, with a constant that depends on the method. In the resolved range a fixed mode converges at the rate $h^2$ for the Galerkin method (theoretical bound $h^p$), $h$ for the nonconsistent method for both $p$, as the $p$-independent resolved order of \cref{thm:supg} predicts, and $h^{p+1}$ for the CIP method; for the stabilised methods this is half an order better than $\sigma-\kappa$. \emph{Dependence on $k$.} In the resolved range the bounds grow like $|k|^\sigma$. On the short window $[1.9,2]$ the coefficients of the stabilised methods grow at exactly this rate (\cref{tab:rates1d}). On $[0,2]$, the window of \cref{fig:burgers_spectra}, the growth is slower, since the residual is concentrated at the shock, whose position moves, so that the phases of the coefficients partly cancel in the time integral: it is close to $|k|^\sigma$ for the CIP method (right column), while the spectrum of the nonconsistent SUPG method (middle column) is nearly flat in the resolved range and lies far above that of the CIP method there, since its stabilisation pollutes all resolved frequencies. The Galerkin coefficients (left column) grow more slowly than $|k|^\sigma$ on both windows.

\emph{Both together.} By \eqref{eq:type}, $h^{\kappa}\rho_k/\Ms_\sigma(h|k|)$, and hence $h^{\kappa}|\Rh(1;\phi_k)|/\Ms_\sigma(h|k|)$, is bounded uniformly in $h$ and $k$; the maximum $C_{\rm obs}$ of the latter over $k$ is given in the last column of \cref{tab:rates1d}. For the Galerkin method $C_{\rm obs}\approx1$ on all meshes: the bound is attained. For the stabilised methods $C_{\rm obs}$ decreases like $h^{1/2}$, and with $\kappa=0$ it is nearly constant ($0.4$ to $1.5$ for $N=32,\dots,512$): the observed spectra are of type $(0,\sigma)$. For the CIP method, with $\sigma=p+1$, this is the benchmark of \cref{sec:burgers}: in this example the CIP method has the optimal residual spectrum, and for $p=2$ its tails, the maxima over $h|k|\ge1$, do not even grow under refinement. For CIP with $p=2$, $\norm{(\mathrm{Id}-\Pi_1^X)u_h}^2\le0.44h$ at the final time, as expected for a shock, consistent with the condition of \cref{prop:cip_highp}.

\begin{table}[!htb]
\centering
\caption{Burgers equation: observed exponents of the residual spectrum and those of the bounds (in parentheses). Peak $\max_{h|k|\ge1}|\Rh(1;\phi_k)|\sim h^{-b}$ and resolved mode $|\Rh(1;\phi_2)|\sim h^{a}$, least-squares fits over $N=64,\dots,512$; growth $|\Rh(\psi;\phi_k)|\sim|k|^{s}$ for $1\le k\le8$, $N=512$, with $\psi\equiv1$ on $[0,2]$ and on $[1.9,2]$; $C_{\rm obs}=\max_k h^{\kappa}|\Rh(1;\phi_k)|/\Ms_\sigma(h|k|)$ with the $(\kappa,\sigma)$ of the bound, for $N=32$ and $N=512$.}\label{tab:rates1d}
\small\setlength{\tabcolsep}{4.5pt}
\begin{tabular}{lcccccc}
\toprule
method & $p$ & $b$ \ ($\kappa$) & $a$ \ ($\sigma-\kappa$) & $s$ on $[0,2]$ & $s$ on $[1.9,2]$ \ ($\sigma$) & $C_{\rm obs}$, $N=32\to512$\\
\midrule
Galerkin & 1 & $0.96$ \ ($1$) & $2.01$ \ ($1$) & $1.57$ & $1.61$ \ ($2$) & $1.09\to0.90$\\
nonconsistent SUPG & 1 & $-0.07$ \ ($\frac12$) & $0.98$ \ ($\frac12$) & $0.15$ & $1.01$ \ ($1$) & $0.34\to0.09$\\
CIP & 1 & $-0.05$ \ ($\frac12$) & $2.01$ \ ($\frac32$) & $1.73$ & $2.00$ \ ($2$) & $0.40\to0.08$\\
\midrule
Galerkin & 2 & $0.97$ \ ($1$) & $1.99$ \ ($2$) & $2.05$ & $2.17$ \ ($3$) & $1.34\to1.24$\\
nonconsistent SUPG & 2 & $0.00$ \ ($\frac12$) & $1.00$ \ ($\frac12$) & $0.15$ & $1.02$ \ ($1$) & $0.18\to0.05$\\
CIP & 2 & $0.00$ \ ($\frac12$) & $2.98$ \ ($\frac52$) & $2.73$ & $2.99$ \ ($3$) & $0.62\to0.16$\\
\bottomrule
\end{tabular}
\end{table}

\subsection{An under-resolved double shear layer}\label{sec:num2d}

We consider the thin double shear layer of \cite{bell1989second,minion1997performance} on the unit torus $(0,1)^2$; the results of \cref{sec:framework,sec:galerkin,sec:stabilised} transfer by scaling, with $h|\bk|$ replaced by $2\pi h|\bk|$. The initial velocity is
\[
u_1(x,y,0)=\begin{cases}\tanh(\varrho(y-\frac14)), & y\le\frac12,\\ \tanh(\varrho(\frac34-y)), & y>\frac12,\end{cases}\qquad u_2(x,y,0)=\delta\sin(2\pi x),
\]
for $\bu=(u_1,u_2)$ and $(x,y)\in(0,1)^2$, with $\varrho=80$, $\delta=0.05$, $\bff=0$, and $\mu=10^{-4}$ \cite{minion1997performance} or $\mu=0$. The layers roll up into two vortices connected by braids which are stretched and become thinner in time; under-resolved computations are known to produce spurious secondary vortices in the braids \cite{minion1997performance}.

To compare the methods on the same velocity space we use quadratic velocities throughout: the Taylor--Hood Galerkin method, its grad--div stabilised version ($\gamma=0.2$), and the nonconsistent streamline-diffusion and CIP methods with $\mathbb P_2/\mathbb P_2$ elements, on structured meshes of $N\in\{16,32,64,128\}$ squares per direction, each divided into two triangles, so that $h=1/N$ (the side of the squares; the element diameter is $\sqrt2h$, which only changes constants), the node spacing is $h_v=h/2$ and the number of velocity nodes per direction is $n_v=2N\in\{32,\dots,256\}$. The initial layers have width $1/80$, that is $0.4$ to $3.2$ node spacings. We take $U=1$ and, in \eqref{eq:supg} with $h$ replaced by the node spacing, $\tau=c_\tau h_v/U$ and $\tau_d=c_dh_vU$ with $c_\tau=c_d=\frac12$, and for CIP $\gamma_u=\gamma_p=0.01$, $\varepsilon=0$ in \eqref{eq:cip_weight} and $\tau_d$ as for the streamline-diffusion method; since $Uh/\mu\ge78$, the pressure weight in \eqref{eq:cip} is $U^{-1}$ in all computations. In time we use the Crank--Nicolson method with extrapolated convecting velocity $\bw^n=\frac32\bu_h^n-\frac12\bu_h^{n-1}$, where $\bu_h^n$ is the approximation at time $n\Delta t$, in $c(\bw^n;\cdot,\cdot)$ and $S_h$, and $\Delta t=\frac12h_v$; with \eqref{eq:skew} this scheme is energy stable, and for the Galerkin method with $\mu=0$ it conserves the energy exactly. The initial data are obtained by a discrete Leray projection: $\bu_{0,h}\in\Vh$ with $(\bu_{0,h},\bv_h)_\Omega-(p,\nabla\cdot\bv_h)_\Omega=(\bu_0,\bv_h)_\Omega$ and $(\nabla\cdot\bu_{0,h},q_h)_\Omega+s(p,q_h)=0$ for all $(\bv_h,q_h)\in\Vh\times\Qh$, where $s$ is the pressure--pressure part of the stabilisation of the method ($s=0$ for the Taylor--Hood methods). The reference solution is computed with a dealiased pseudo-spectral method on a $512^2$ grid; it differs from the $256^2$ solution by less than $4\cdot10^{-7}$ in all modes with $|\bk|\le16$.

\begin{figure}[!htb]
\centering
\includegraphics[width=\textwidth]{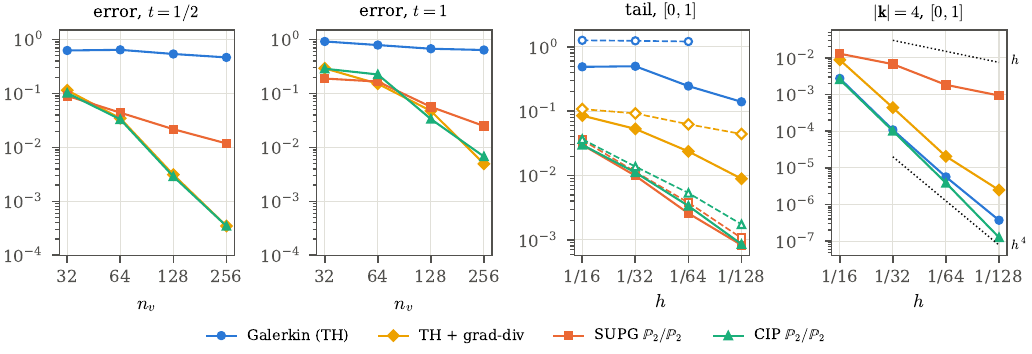}
\caption{Double shear layer, from left to right: error in the modes $0<|\bk|\le8$ at $T=\frac12$ and $T=1$ against $n_v$; maximum of the momentum residual over $2\pi h|\bk|\ge1$ (tail) for $\mu=10^{-4}$ (solid) and $\mu=0$ (dashed; Galerkin method $N\le64$), and momentum residual at $|\bk|=4$ (dotted: slopes $h$ and $h^4$), both on $[0,1]$, against $h$.}\label{fig:shear_errors}
\end{figure}

\begin{figure}[!htb]
\centering
\includegraphics[width=\textwidth]{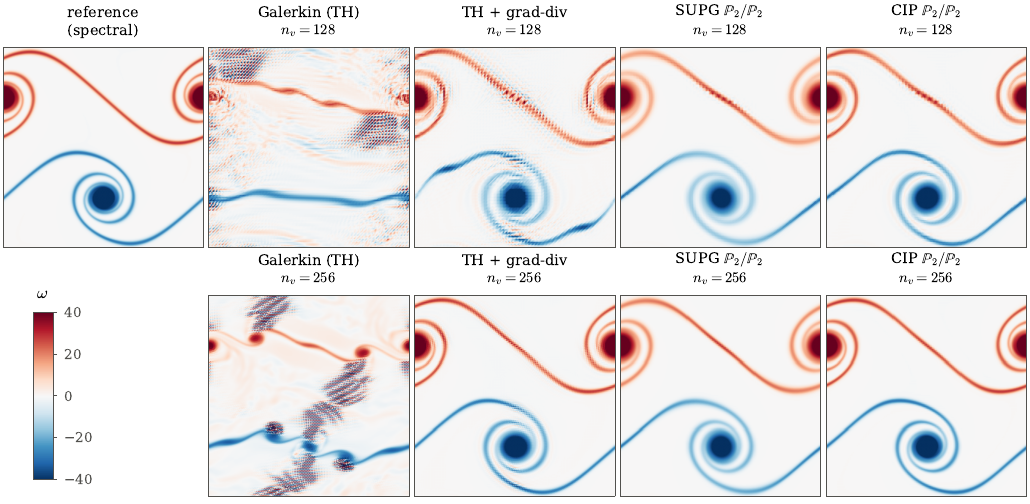}
\caption{Double shear layer, $\mu=10^{-4}$: vorticity $\omega=\partial_xu_2-\partial_yu_1$ at $t=1$ for the reference solution and the four methods with $n_v=128$ (top, $N=64$) and $n_v=256$ (bottom, $N=128$) velocity nodes per direction.}\label{fig:shear_vorticity}
\end{figure}

\begin{figure}[!htb]
\centering
\includegraphics[width=\textwidth]{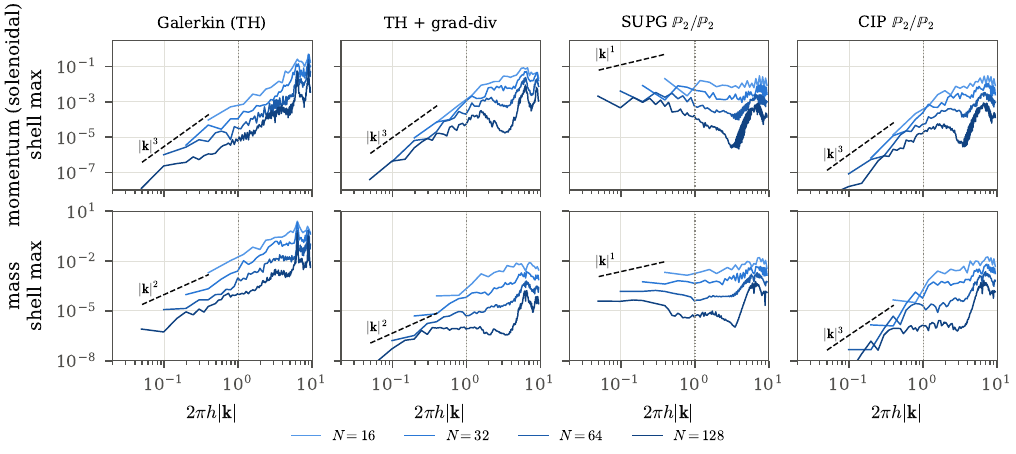}
\caption{Double shear layer, $\mu=10^{-4}$: shell maxima of the time-weak residual spectra ($\psi\equiv1$ on $[0,1]$) against $2\pi h|\bk|$; top: solenoidal momentum coefficients, bottom: mass coefficients. Dotted line: $2\pi h|\bk|=1$; dashed: the growth $|\bk|^\sigma$ of the bounds (\cref{tab:summary}).}\label{fig:shear_spectra}
\end{figure}

We consider the error in the Fourier modes $0<|\bk|\le8$ of the velocity at $T=\frac12$ and $T=1$, $(\sum_{0<|\bk|\le8}|\hat\be_{\bk}(T)|^2)^{1/2}$ with $\hat\bv_{\bk}:=\int_{(0,1)^2}\bv\,e^{-2\pi i(k_1x+k_2y)}\dd x\dd y$ (the full vector coefficients; their longitudinal part, which the representation with solenoidal terminal data of \cref{sec:framework} does not cover, is at most $6.2\cdot10^{-4}$ times the error at $T=1$ for $n_v=128$ and $256$), and the corresponding residual spectra with $\psi\equiv1$ on $[0,T]$: for each $\bk$ with $|k_1|,|k_2|\le\frac32N$ we compute the solenoidal momentum coefficient ($\ba=(-k_2,k_1)/|\bk|$) and the mass coefficient, and report shell maxima over $|\,|\bk|-n|<\frac12$, the shell maximum for $n=4$ (``$|\bk|=4$'', resolved for $N\ge32$) and the maximum over $2\pi h|\bk|\ge1$ (tail). At $T=\frac12$ the vortices have formed and the braids have about the width of the initial layers; between $T=\frac12$ and $T=1$ they are stretched further and become unstable, so that the low modes at $T=1$ are much more sensitive to perturbations, that is, the adjoint weights $\Phi_{\bk}$ in \eqref{eq:spectral_bound} are larger. The results are collected in \cref{tab:shear} and \cref{fig:shear_vorticity,fig:shear_spectra,fig:shear_errors}.

\begin{table}[!htb]
\centering
\footnotesize
\caption{Double shear layer: error in the modes $0<|\bk|\le8$ at time $T$, momentum residual at $|\bk|=4$ and its maximum over $2\pi h|\bk|\ge1$ (tail), with $\psi\equiv1$ on $[0,T]$, for $\mu=10^{-4}$; last column: tail for $\mu=0$ (--: not computed). The error is $(\sum_{0<|\bk|\le8}|\hat\be_{\bk}(T)|^2)^{1/2}$.}\label{tab:shear}
\setlength{\tabcolsep}{3pt}
\begin{tabular}{lrrccccccc}
\toprule
 & & & \multicolumn{3}{c}{$T=\frac12$} & \multicolumn{4}{c}{$T=1$}\\
\cmidrule(lr){4-6}\cmidrule(lr){7-10}
method & $N$ & $n_v$ & error & $|\bk|=4$ & tail & error & $|\bk|=4$ & tail & tail, $\mu=0$\\
\midrule
Galerkin TH & 16 & 32 & 0.628 & $1.5\cdot10^{-3}$ & 0.448 & 0.923 & $2.7\cdot10^{-3}$ & 0.489 & 1.263\\
 & 32 & 64 & 0.645 & $1.6\cdot10^{-4}$ & 0.493 & 0.791 & $1.1\cdot10^{-4}$ & 0.499 & 1.241\\
 & 64 & 128 & 0.541 & $4.5\cdot10^{-6}$ & 0.236 & 0.676 & $5.6\cdot10^{-6}$ & 0.246 & 1.211\\
 & 128 & 256 & 0.467 & $2.3\cdot10^{-7}$ & 0.130 & 0.641 & $3.7\cdot10^{-7}$ & 0.140 & --\\
\midrule
TH + grad--div & 16 & 32 & 0.115 & $4.4\cdot10^{-3}$ & 0.0515 & 0.293 & $8.8\cdot10^{-3}$ & 0.0853 & 0.108\\
 & 32 & 64 & 0.0347 & $2.6\cdot10^{-4}$ & 0.0291 & 0.152 & $4.4\cdot10^{-4}$ & 0.0531 & 0.0923\\
 & 64 & 128 & $3.1\cdot10^{-3}$ & $1.0\cdot10^{-5}$ & 0.0108 & 0.0476 & $2.1\cdot10^{-5}$ & 0.0237 & 0.0626\\
 & 128 & 256 & $3.5\cdot10^{-4}$ & $7.8\cdot10^{-7}$ & $4.4\cdot10^{-3}$ & $5.0\cdot10^{-3}$ & $2.5\cdot10^{-6}$ & $8.9\cdot10^{-3}$ & 0.0443\\
\midrule
SUPG, $\mathbb P_2/\mathbb P_2$ & 16 & 32 & 0.0905 & $2.0\cdot10^{-3}$ & 0.0223 & 0.189 & 0.0129 & 0.0300 & 0.0357\\
 & 32 & 64 & 0.0441 & $1.8\cdot10^{-3}$ & $9.6\cdot10^{-3}$ & 0.166 & $6.7\cdot10^{-3}$ & $9.9\cdot10^{-3}$ & 0.0114\\
 & 64 & 128 & 0.0218 & $8.6\cdot10^{-4}$ & $3.0\cdot10^{-3}$ & 0.0565 & $1.8\cdot10^{-3}$ & $2.6\cdot10^{-3}$ & $3.8\cdot10^{-3}$\\
 & 128 & 256 & 0.0119 & $4.0\cdot10^{-4}$ & $8.6\cdot10^{-4}$ & 0.0249 & $9.4\cdot10^{-4}$ & $8.3\cdot10^{-4}$ & $1.1\cdot10^{-3}$\\
\midrule
CIP, $\mathbb P_2/\mathbb P_2$ & 16 & 32 & 0.102 & $9.1\cdot10^{-4}$ & 0.0246 & 0.288 & $2.6\cdot10^{-3}$ & 0.0300 & 0.0367\\
 & 32 & 64 & 0.0331 & $5.0\cdot10^{-5}$ & $9.9\cdot10^{-3}$ & 0.225 & $1.0\cdot10^{-4}$ & 0.0112 & 0.0138\\
 & 64 & 128 & $2.9\cdot10^{-3}$ & $1.7\cdot10^{-6}$ & $3.3\cdot10^{-3}$ & 0.0336 & $3.9\cdot10^{-6}$ & $3.4\cdot10^{-3}$ & $5.4\cdot10^{-3}$\\
 & 128 & 256 & $3.5\cdot10^{-4}$ & $5.0\cdot10^{-8}$ & $8.7\cdot10^{-4}$ & $6.8\cdot10^{-3}$ & $1.3\cdot10^{-7}$ & $8.6\cdot10^{-4}$ & $1.7\cdot10^{-3}$\\
\bottomrule
\end{tabular}
\end{table}

\emph{Stable quantity, $T=\frac12$.} Here the errors behave as predicted by the residual spectra. The Galerkin method does not converge on these meshes (errors between $0.65$ and $0.47$), although its resolved residual is small, because its unresolved residuals are large (momentum tails $0.13$ to $0.49$, mass tails up to $1.9$); for $\mu=0$ they do not decrease under refinement, consistent with the plateau $h^{-1}$ of \cref{thm:th}. The streamline-diffusion method converges at first order, as does its resolved residual, in agreement with the $p$-independent resolved order of \cref{thm:supg}. The CIP and grad--div methods converge fast, both with error $3.5\cdot10^{-4}$ for $n_v=256$. The CIP method has the best residual spectrum in the sense of \cref{tab:summary}: on all meshes the smallest resolved momentum residual of all methods, decreasing by factors of $18$ to $33$ per refinement (its mass residual is larger than that of the grad--div method on most shells), and, together with the streamline-diffusion method, the smallest tails, which are $2$ to $5$ times smaller than those of the grad--div method and $18$ to $150$ times smaller than those of the Galerkin method. This ordering of the tails is that of the tail exponents in \cref{tab:summary}: $\kappa=\frac12$ for the streamline-diffusion and CIP methods and $\kappa=1$ for the Taylor--Hood methods, for which grad--div stabilisation improves the bound of the mass residual only. It is also reflected in the dependence on the viscosity: the tails of the streamline-diffusion and CIP methods for $\mu=0$ exceed those for $\mu=10^{-4}$ by at most a factor $2$, so that the dissipation of the unresolved scales is provided by the stabilisation, whereas those of the grad--div and Galerkin methods increase by factors up to $5.0$ (on $[0,1]$). For $\mu=0$ the tail of the grad--div method decreases only like $h^{0.5}$, against $h^{1.6}$ for the CIP method, but its vorticity stays bounded and its divergence decreases, consistent with the conditions of \cref{rem:gd_curl}. A posteriori, at $t=\frac12$ and $t=1$, on all meshes and for both viscosities, the CIP solutions satisfy the convective part of (S1) with $C_S=1.4$ to $1.7$ (snapshot values, consistent with the time-integrated condition of \cref{prop:cip_highp}), the elementwise part $(\mathrm{Id}-\pi^X)(\bu_h\cdot\nabla)\bu_h$ contributing at most $10\%$, and $h^{-1/2}\norm{(\mathrm{Id}-\Pi_1^X)\bu_h}^2\le4\cdot10^{-3}\norm{\bu_h}^2$, so that the additional term of \cref{prop:cip_highp} is negligible; both decrease under refinement.

\emph{Orders in $h$ and $k$.} We test the bounds as in \cref{sec:num1d}, for the momentum residual on $[0,1]$ with $\mu=10^{-4}$ (\cref{fig:shear_spectra}). In the resolved range $2\pi h|\bk|\le\frac12$ the shell maxima grow like $|\bk|^{2.7}$ to $|\bk|^{3.0}$ for the grad--div and CIP methods, as the factor $|\bk|^{p+1}$ of the bounds, and more slowly for the Galerkin method (about $|\bk|^{1.7}$) and the streamline-diffusion method, whose coefficients hardly depend on $|\bk|$. In $h$ they decrease much faster than the bounds: at fixed $|\bk|=4$ at the rates $3.7$ to $4.8$ for the Taylor--Hood and CIP methods (\cref{fig:shear_errors}, right) and $0.7$ to $1.4$ for the streamline-diffusion method, about half an order better than $\sigma-\kappa=\frac12$, as in \cref{sec:num1d}; the plateaus decrease like $h^{1.7}$ for the streamline-diffusion and CIP methods, $h^{1.1}$ for the grad--div method and $h^{0.6}$ for the Galerkin method, whose plateau stays at about $1.2$ for $\mu=0$. Accordingly, $h^{\kappa}$ times the computed shell maxima divided by $\Ms_\sigma(2\pi h|\bk|)$ decreases under refinement for all methods. The bounds hold, but are attained only in a regime that remains under-resolved, such as the shock of \cref{sec:num1d}; here the layers become resolved on the finer meshes, and the residual decays at the faster rates allowed by the smoothness of the solution, which bounds based only on the energy estimate cannot see.

\emph{Time discretisation.} The resolved residual coefficients quoted above, and their rates, are those of the algebraic residual \eqref{eq:discrete_residual}. 
In the errors, the time discretisation becomes visible on the finest meshes. At $T=1$ the time error is at most $4\%$ of the errors, and dividing $\Delta t$ by four at $n_v=64$ changes the CIP errors by less than $2\%$.

\emph{Sensitive quantity, $T=1$.} On $[0,1]$ the residual spectra are ordered as on $[0,\frac12]$, but on the meshes with $n_v\le64$ the errors are not: those of the grad--div and CIP methods, which dissipate least, are dominated by fully developed spurious secondary vortices in the braids, due to under-resolution. For $n_v=128$ the grad--div method shows grid-scale wiggles along the braids instead (\cref{fig:shear_vorticity}, top, middle), these wiggles are less accentuated in the CIP and streamline-diffusion methods. The more dissipative streamline-diffusion method also smears the braids, reflecting its slow convergence of the resolved scales of the residual. On the finest mesh, $n_v=256$, the vortices and wiggles have disappeared for the grad--div, CIP and streamline-diffusion methods (\cref{fig:shear_vorticity}, bottom), and the errors, $5.0\cdot10^{-3}$ (grad--div), $6.8\cdot10^{-3}$ (CIP) and $0.025$ (streamline diffusion), fall into the same classes as at $T=\frac12$, while the Galerkin method hardly improves ($0.68$ and $0.64$ for $n_v=128$ and $256$). This is consistent with \cref{thm:abstract}: the adjoint weights for $T=1$ are large in the unstable braids, so that the residual must be small up to higher wave numbers before the error becomes small. A smaller residual spectrum gives a smaller bound; it orders the errors only when the quantity is stable enough for the bound to be informative.

\emph{Divergence consistency term.} At $t=1$, $\mu=10^{-4}$ and $n_v=256$, the largest coefficient of $\frac12(\nabla\cdot\bu_h)\bu_h$, the density of $\epsilon_h$ in \eqref{eq:Kh}, for the solenoidal modes $\ba\phi_{\bk}$, $0<|\bk|\le8$, is $2.6\cdot10^{-3}$ (Galerkin), $1.9\cdot10^{-6}$ (grad--div), $4.5\cdot10^{-5}$ (streamline diffusion) and $7.3\cdot10^{-7}$ (CIP), between $870$ and $7.7\cdot10^4$ times below the bound $\frac12\norm{\nabla\cdot\bu_h}\norm{\bu_h}$ that \eqref{eq:Kh_general} uses for each coefficient; for CIP it decreases from $2.8\cdot10^{-3}$ to $7.3\cdot10^{-7}$ for $n_v=32,\dots,256$, while this bound only decreases from $0.18$ to $0.056$.

\subsection{Forced turbulence far from the viscous range}\label{sec:forced}
\begin{figure}[!htb]
\centering
\includegraphics[width=\textwidth]{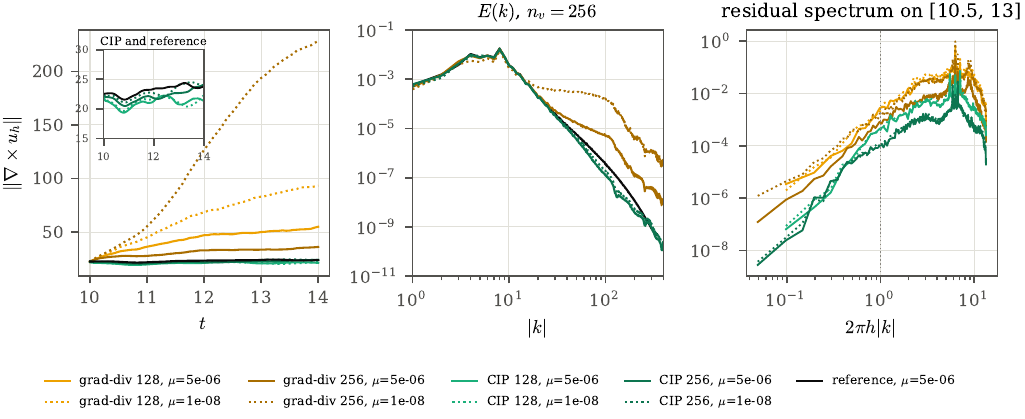}
\caption{Forced turbulence, grad--div and CIP methods with $n_v=128$ and $256$, $\mu=5\cdot10^{-6}$ (solid) and $\mu=10^{-8}$ (dotted), and the pseudo-spectral reference ($\mu=5\cdot10^{-6}$, black). Left: $\norm{\nabla\times\bu_h}$ against $t$ (inset: CIP and reference). Middle: energy spectrum $E(k)=\frac12\sum_{|\,|\bk|-k|<1/2}|\hat\bu_{\bk}|^2$, averaged over $[10.5,13]$, for $n_v=256$. Right: shell maxima of the momentum residual on $[10.5,13]$ ($\psi\equiv1$) against $2\pi h|\bk|$.}\label{fig:forced}
\end{figure}

The double shear layer is not forced, so that its vorticity is bounded by that of the initial data, which favours the grad--div method (\cref{rem:gd_curl}). We have therefore also computed forced two-dimensional turbulence, in which enstrophy is transferred to the small scales at a constant mean rate: $\bff=(\sin(16\pi y),0)$, a drag $-\alpha\Pi_{<8}\bu$ with $\alpha=1$, where $\Pi_{<8}$ is the $L^2$ projection onto the Fourier modes $0<|\bk|<8$ (it arrests the inverse energy cascade and leaves the forced modes and all smaller scales undamped), and $\mu=5\cdot10^{-6}$, with the discretisations of \cref{sec:num2d}. All computations start at $t=10$ from a pseudo-spectral solution in statistical equilibrium; the pseudo-spectral solution on a $1024^2$ grid, whose viscous cut-off is at $|\bk|\approx160$, serves as reference. Since $\nabla\times$ commutes with $\Pi_{<8}$ and $((\bu\cdot\nabla)\omega,\omega)_\Omega=0$, the vorticity $\omega=\nabla\times\bu$ of the exact solution satisfies $\frac12\frac{\dd}{\dd t}\norm{\omega}^2\le(\nabla\times\bff,\omega)_\Omega$ for every $\mu\ge0$, hence
\[
\norm{\omega(t)}\le\norm{\omega(10)}+(t-10)\norm{\nabla\times\bff},\qquad\norm{\nabla\times\bff}=8\sqrt2\,\pi\approx35.5,\quad\norm{\omega(10)}\approx22.5.
\]
Averaged over $t\in[10.5,13]$, $\norm{\nabla\times\bu_h}$ of the grad--div method is $1.87$ and $1.35$ times the reference value $22.9$ for $n_v=128$ and $256$, and the tail of its momentum residual on $[10.5,13]$ ($\psi\equiv1$, forcing and drag included) is $0.51$ and $0.21$; for the CIP method these numbers are $0.92$ and $0.95$, and $0.072$ and $0.028$. The improvement of the grad--div method on the finer mesh is due to the viscosity: its grid scale $|\bk|=128$ is close to the viscous cut-off, where $\mu(2\pi|\bk|)^2\approx3$, so that this computation is an under-resolved direct numerical simulation rather than a large eddy simulation. With $\mu=10^{-8}$, for which $\mu(2\pi|\bk|)^2\approx6.5\cdot10^{-3}$ at the grid scale, $\norm{\nabla\times\bu_h}$ of the grad--div method becomes $2.7$ and $4.8$ times the same reference value (that of $\mu=5\cdot10^{-6}$; there is no reference solution for $\mu=10^{-8}$) and the tail $0.76$ and $0.99$: both grow under refinement. For $n_v=256$, $\norm{\nabla\times\bu_h}$ increases from $22.5$ at $t=10$ to $228$ at $t=14$, where the bound above is $165$; it exceeds this bound from $t\approx11.2$ on, so that the growth is not physical. The CIP method is insensitive to the change of viscosity: $0.93$ and $0.97$ times the reference value, with tails $0.073$ and $0.028$, for $n_v=128$ and $256$; for $n_v=256$, $\norm{\nabla\times\bu_h}$ stays below $24.5$ up to $t=14$. The time evolution of the vorticity and the spectra are reported in \Cref{fig:forced}. In the enstrophy range $8<|\bk|\lesssim160$ the reference spectrum decays like $|\bk|^{-4}$ (least-squares exponents $-3.8$ to $-4.0$ for $15\le|\bk|\le120$), steeper than the prediction $|\bk|^{-3}$ of \cite{kraichnan1967inertial,leith1968diffusion,batchelor1969computation} and close to the $|\bk|^{-4}$ that Saffman~\cite{saffman1971spectrum} derived for vorticity distributions with discontinuities; such spectra are typical of forcing that is constant in time, which produces strong long-lived vortices, at moderate resolution \cite[Sect.~4.3]{boffetta2012two}. The spectrum of the CIP method follows it, whereas that of the grad--div method flattens above $|\bk|\approx20$. These computations indicate that, as long as the mesh is far from the viscous range, the computed enstrophy and the residual tail of the grad--div method grow under refinement, so that the first condition of \cref{rem:gd_curl} can fail also in two dimensions, whereas the gradient jumps of the CIP method control the grid-scale vorticity.

\section{Concluding remarks}\label{sec:conclusions}

We have proposed to assess computational methods for large eddy simulation through the spectrum of their residual. The criterion separates the method-dependent part of a dual weighted error representation, the residual, from the problem-dependent part, the adjoint solution; the residual spectrum is bounded using only the energy stability of the method, without regularity of the exact solution and uniformly in the viscosity, and it can be computed a posteriori. Two exponents characterise a method: the resolved order $\sigma-\kappa$ and the tail exponent $\kappa$, which determines how stable a quantity of interest must be for its approximation to converge. Of the methods considered, the continuous interior penalty method with grad--div stabilisation combines the highest resolved order, $p+\frac12$, with the smallest tail exponent, $\frac12$ (for $p\ge2$ under an a posteriori condition); the Galerkin method has tail exponent $1$, grad--div stabilisation improves its resolved order but not its tail, and the nonconsistent streamline-diffusion method has the small tail exponent but a resolved order limited by its lack of consistency.

The computations confirm these predictions.

\Cref{tab:ranking} condenses the analysis and the computations into an informal ranking of the methods for implicit large eddy simulation, without regard to cost. In our computations the control of the grid scales by the CIP method did not depend on the viscosity. The order of the second and third places depends on the regime: the grad--div method is far more accurate than the streamline-diffusion method when the flow is decaying or the mesh approaches the viscous range, but its control of the grid-scale vorticity is not guaranteed under sustained forcing far from the viscous range. For a sensitive quantity no method is reliable until the mesh resolves the scales to which the quantity is sensitive.

\begin{table}[ht]
\centering
\small
\caption{Informal ranking of the methods for implicit large eddy simulation, without regard to cost. Type $(\kappa,\sigma)$ from \cref{tab:summary} ($^*$for $p\ge2$ under the condition of \cref{prop:cip_highp}; $^\dagger(\frac12,p+1)$ a posteriori under the conditions of \cref{rem:gd_curl}). Error: double shear layer, modes $0<|\bk|\le8$, $T=\frac12$, $n_v=256$, $\mu=10^{-4}$ (\cref{tab:shear}). Forced: forced turbulence of \cref{sec:forced}, $n_v=256$, time average of $\norm{\nabla\times\bu_h}$ relative to the reference, and residual tail on $[10.5,13]$; $\mu=5\cdot10^{-6}$ unless stated; --: not computed.}\label{tab:ranking}
\begin{tabular}{@{}c>{\raggedright\arraybackslash}p{2.3cm}cc>{\raggedright\arraybackslash}p{2.6cm}>{\raggedright\arraybackslash}p{4.9cm}@{}}
\toprule
 & method & $(\kappa,\sigma)$ & error & forced & assessment\\
\midrule
1 & CIP, bulk grad--div & $(\frac12,p+1)^*$ & $3.5\cdot10^{-4}$ & $0.95$, $0.028$ & best residual spectrum; grid-scale vorticity and tail controlled in all computations, independently of the viscosity\\
\addlinespace
2 & TH + grad--div & $(1,p+1)^\dagger$ & $3.5\cdot10^{-4}$ & $1.35$, $0.21$; for $\mu=10^{-8}$: $4.8$, $0.99$ & accurate for decaying flows and near the viscous range; under forcing far from it, unphysical growth of the grid-scale vorticity and of the tail under refinement\\
\addlinespace
3 & nonconsistent SUPG & $(\frac12,1)$ & $0.012$ & -- & small tails and robust, but first order in the resolved range\\
\addlinespace
4 & Galerkin TH & $(1,p)$ (mass) & $0.47$ & -- & no control of the grid scales; does not converge on the meshes considered\\
\bottomrule
\end{tabular}
\end{table}

We expect results similar to those of the CIP method for all methods with the classical $h^{p+1/2}$ estimate for smooth solutions, such as local projection and orthogonal subscale stabilisation, and for $H(\mathrm{div})$-conforming and discontinuous Galerkin methods with upwind fluxes; Petrov--Galerkin methods, whose test functions depend on the discrete velocity, are harder to interpret in this framework (the nonconsistent method of \cref{sec:supg} is not of this kind: its stabilisation is symmetric). Finally, closing the argument for turbulent flows requires information on the spectrum of the adjoint solution; \cref{thm:abstract,thm:negative} show precisely how it enters, and sharper, possibly statistical, estimates of adjoint spectra for model flows would allow the argument to be closed in specific situations. A companion paper will consider finite element methods with a divergence-free discrete velocity subspace, such as the Scott--Vogelius element and $H(\mathrm{div})$-conforming methods, for which the divergence consistency term vanishes identically.

\paragraph{Code availability} The code reproducing all numerical examples in this paper is available from the authors upon request.

\paragraph{Acknowledgement.} This work was funded by the EPSRC grants EP/T033126/1 and EP/X042650/1. During the preparation of this manuscript, the authors used Anthropic Claude Opus 5 to assist with drafting and revising text, improving mathematical exposition, checking notation and internal consistency, editing LaTeX, and developing computational code. 
All AI-assisted material, including mathematical statements and computational code, was critically reviewed and validated by the authors, who take full responsibility for the accuracy, originality, and integrity of the manuscript.

\bibliographystyle{abbrvnat}
\bibliography{refs}

\end{document}